\documentclass[10pt]{article}
\usepackage{ amssymb, amsmath}
\usepackage{amsfonts}
\usepackage{amsthm}
\usepackage{enumerate}
\newtheorem{theorem}{Theorem}[section]

\newtheorem{proposition}[theorem]{Proposition}

\newtheorem{definition}[theorem]{Definition}
\newtheorem{remark}[theorem]{Remark}
\newcommand{\bysame}{\mbox{\rule{3em}{.4pt}}\,}

\def\dim{\textrm{ dim }}
\def\Ker{\textrm { Ker }}
\def\lra{\longrightarrow}
\def\p{\partial}
\def\rank{\textrm{ rank }}
\def\R{\mathbb R}
\def\so{{\cal S}_1}
\begin{document}
\topmargin=-.3in

\oddsidemargin=.15in

\evensidemargin=.15in

\title{An FIO calculus for marine seismic imaging: folds and cross caps}

\author{ Raluca Felea  and  Allan Greenleaf
\thanks{The second author was partially supported by NSF grants
DMS-0138167 and DMS-0551894.} }
\date{}

\maketitle

\begin{abstract}
We consider a linearized inverse problem arising in offshore seismic imaging.
Following Nolan and Symes\cite{nosy}, one wishes to determine a
singular perturbation
of a smooth background soundspeed in the Earth from measurements made
at the surface
resulting from various seismic experiments; the overdetermined data
set considered
here corresponds to marine seismic exploration. In the presence of only
fold caustics for the background, we identify the  geometry of the
canonical relation
underlying the linearized forward scattering operator $F$, which is a Fourier
integral operator. We then establish a composition calculus for general FIOs
associated with  similar canonical relations, which we call
\emph{folded cross caps},
sufficient for identifying the normal operator $F^*F$. In contrast to
the case of a
single source experiment, treated by Nolan\cite{no} and Felea\cite{fel}, the
resulting artifact is
$\frac12$ order smoother than the main pseudodifferential part of $F^*F$.

\end{abstract}

\section{Introduction}

\par This article deals with a linearized inverse scattering problem
considered by  Nolan and  Symes [23]. Acoustic waves are generated
at the surface of the earth, scatter off heterogeneities in the
subsurface and return to the surface. The full inverse problem
would use the pressure field at the surface to reconstruct an
image of the subsurface. We instead consider the linearized
operator $F$ which maps singular perturbations of a smooth
background  sound speed in the subsurface, assumed known, to
perturbations of the resulting pressure field at the surface. The
goal is to left-invert $F$; standard techniques suggest  studying
left invertibility of the normal operator $N=F^*F$. To start, we
make two assumptions: $(i)$ no single ray connects a source to a
receiver; and $(ii)$ no ray originating in the subsurface grazes
the surface. Under these assumptions, in the case of a {\it
single} source and receivers ranging over an open subset of the
surface, $\{x_3=0\}$, Rakesh\cite{ra} showed that $F$ is a Fourier integral
operator (FIO).  Beylkin\cite{be} showed that if caustics do not occur
for the background soundspeed, $F^*F$ is a pseudodifferential
operator $ (\Psi DO)$.

\par For more general data acquisition geometries,  the canonical
relation of $F$ depends on the sets of sources  and receivers.
Nolan and Symes\cite{nosy} proved that, if both sources and receivers
vary over open and bounded subsets $\Sigma_r$ and $\Sigma_s $ of the
surface, then under the traveltime injectivity condition (TIC),
generalizing the no-caustic assumption, $F^*F$ is still a
$\Psi DO$. The same result was stated by ten Kroode, Smit and
Verdel\cite{ksv} and their proof was completed by Stolk\cite{st} who also
relaxed the TIC condition in low dimensions.

\par For applications in three spatial variables, an important problem is to
understand the nature of $F$ and $F^*F$ for the {\it marine} data acquisition
geometry \cite{nosy}, where measurements are made on the codimension
one submanifold
$\Sigma_{r,s}=\{ (r_1, r_2; s_1, s_2) \in \Sigma_r \times
\Sigma_s: s_2=r_2 \}$. This arises as follows: a seismic
vessel trails behind it both an acoustic source and recording
instruments. The point
source consists of an airgun which sends acoustic waves through the
ocean to the
subsurface. Reflections occur when the sound waves encounter
singularities  in the material of the subsurface. The reflected
rays are received by a linear array of hydrophones towed behind the
vessel. The vessel then makes repeated passes along parallel lines
(say, parallel to $x_1$ axis).

\par The purpose of this paper is to  consider  the marine
geometry  under the assumption that only the simplest, most prevalent type of
caustics, namely {\it fold caustics}, occur for the background soundspeed.
Fold caustics are initially defined as follows: A ray departing from
a source $s$ in
the direction $\alpha$
reaches at time $t$ a point denoted $x(t,\alpha)$ in the subsurface. If
there is a
source $s$, such that the spatial projection map $(t,\alpha ) \rightarrow x(t,
\alpha)$ has  a
fold singularity and only singularities of this type, then we say that the
background soundspeed exhibits a fold caustic.  By the
stability of folds, the maps  $(t,\alpha ) \rightarrow x(t,
\alpha)$ also have at most fold singularities for all nearby
sources $s'$. However, it seems that the natural notion of a fold
caustic in the
context of the overdetermined marine data set considered here is the
requirement that
the analogous spatial projection be a {\it submersion with folds}, which is the
simplest singularity in the non-equidimensional setting. This will be
elaborated  upon
in \S2  and \S4.

\par We now introduce the linearized scattering operator $F$
considered in \cite{nosy},\cite{ksv}. The model for the scattered
waves is given by
the wave equation:

\begin{eqnarray}\label{star}
\frac{1}{c^2(x)} \frac{\partial^2p}{\partial
t^2}(x,t)-\triangle p (x,t)&=&\delta(t)\delta(x-s)\\
p(x,t)&=&0,\ \ t < 0,\nonumber
\end{eqnarray}
where $x \in Y=\mathbb R^3_+= \{x \in \R^3, x_3 \geq 0 \}$
represents the Earth, $p(x,t)$
is the pressure field resulting from a pulse at  the source $s$
and $c(x)$ is the velocity field. The linearization consists in
assuming $c$ to be of the form $c=c_0 + \delta c$, where $c_0 $ is
a smooth known background field. The associated  pressure field
$p_0$ is also assumed known. The linearization of (\ref{star}) then
becomes

\begin{eqnarray}\label{starstar}
\frac{1}{c_0^2(x)}\frac{\partial^2 \delta
p}{\partial  t^2}(x,t) -\triangle \delta p(x,t)&=& \frac{2\delta
c(x)}{c_0^3(x) } \frac{\partial^2 p_0}{\partial t^2}\\
   \delta p &=&0,\ \ t < 0,\nonumber
\end{eqnarray}
where $p=p_0 + \delta p$. Now, for a given
data acquisition
submanifold  $\Sigma_{r,s} \subset \partial Y \times \partial Y$
and appropriate time interval $(0,T)$, we define the linearized
scattering operator  $F : \delta c \rightarrow \delta
p|_{\Sigma_{r,s} \times (0,T)}$. The assumption $(ii)$ ensures that
$F$ is an FIO (\cite{ha},\cite{ksv},\cite{ra},\cite{nosy}) and (i)
ensures that the composition $F^*F$ makes sense.

\par  In the case
of the single source model, with only fold caustics appearing,
Nolan\cite{no} showed that $F$ is an FIO associated to a folding
canonical relation in the sense of \cite{meta} (also called a two-sided fold),
and stated that
the Schwartz kernel of the operator $F^*F$ belongs to a class of
distributions associated to two cleanly intersecting Lagrangians
in $(T^*Y \setminus 0) \times (T^*Y \setminus 0)$. This was fully
proved in \cite{fel}. The corresponding canonical relations are the
diagonal $\Delta$ and a  folding canonical relation, different from
the original
one, which lies in $T^*X\times  T^*Y$.

In this article we show that, for the
the marine geometry, the linearization
$F$ is an FIO associated to what we call a {\it folded cross cap }
canonical relation. We then prove  that the Schwartz kernel of $F^*F$
belongs to a class of distributions with a microlocal structure
similar to that in the case of the  single source geometry, but with
the order of the non-pseudodifferential part of $F^*F$ being
$\frac{1}{2}$ lower than in the case of a single source.
\emph{This means that  artifacts arising in  seismic imaging
from the presence
of  fold caustics are $\frac{1}{2}$ derivative smoother for the marine
geometry than for the single source geometry.}

\par  Composition of FIOs under other singular geometries arising in integral
geometry and inverse problems has been previously studied in, e.g.,
\cite{gu-book},\cite{gruh1},\cite{gruh3},\cite{gruh4},\linebreak\cite{no} and
\cite{fel}.

\par The article is organized as follows. In $\S 2$
we review some $C^\infty$ singularity theory and define the
submersion with folds
and  cross cap singularities.
   $\S 3$ is a review
of the distribution classes associated to two cleanly intersecting
Lagrangians, $I^{p,l}(\Lambda_0, \Lambda_1)$ and the operators which
have these as
their Schwartz kernels. In
$\S4$ we show that  submersions with folds and cross caps appear
microlocally in the
marine geometry in the presence of the fold caustics, and we
formulate a general
class of canonical relations exhibiting these singularities.
$\S5$ is dedicated to analyzing a model folded cross cap canonical relation,
$C_0$, in $T^*\R^n \times T^*\R^{n-1}$ ;  we establish  the
composition calculus for $F^*F$, showing  that $F^*F \in I^{p,l}
(\Delta, \tilde{C}_0)$ where $\tilde{C}_0$ is a folding canonical
relation. Finally, \S6 provides the extension of this to the general
class of folded
cross caps. We find a weak normal form for any folded cross cap
canonical relation $C \subset T^*X \times T^*Y$ which allows us to
show that  $F^*F \in I^{p,l} (\Delta, \tilde{C})$, with
$\tilde{C}$ a  folding canonical relation in $T^*Y \times T^*Y$.

\par We would like to thank Cliff Nolan for the helpful discussions, clarifying
\cite{no},
   at the Institute for Mathematics and
its Applications, Minneapolis, in October, \nolinebreak2005.

\section{Fourier integral operators  and singularity classes}

\par Let  $X$  and $Y$ be manifolds and $I^m(X,Y;C)$  denote the class of
$m$-th order Fourier
integral operators (FIOs), $F: {\cal E}'(Y) \to {\cal D}'(X)$, associated
to a canonical relation $C \subset
(T^*X  \setminus 0) \times (T^*Y \setminus 0)$. We will focus on
the composition calculus for two FIOs. Let $C_1 \subset (T^*X
\setminus 0) \times (T^*Y \setminus 0)$ and $C_2 \subset (T^*Y
\setminus 0) \times (T^*Z \setminus 0)$ be two canonical relations
and $F_1 \in I^{m_1}(X,Y;C_{1})$ and $F_2 \in I^{m_2}(Y,Z;C_{2})$.
If $ C_1 \times C_2 $ intersects $T^*X \times \Delta_{T^*Y}\times
T^*Z$ transversally, then H\"{o}rmander\cite{ho} proved that $F_1 \circ
F_2 \in I^{m_1+m_2}(X,Z;C_{1} \circ C_{2})$ where $C_{1} \circ
C_{2}$ is
the composition of $C_1$ and $C_2$ as relations in $T^*X \times T^*Y$
and $T^*Y \times T^*Z$. Duistermaat and Guillemin\cite{dugu} and
Weinstein\cite{we}
extended this calculus to the case of clean
intersection and showed that if $C_1 \times C_2$ and $T^*X \times
\Delta_{T^*Y} \times T^*Z$ intersect
cleanly with
excess $e$ then $A \circ B \in I^{m_1+m_2+e/2}(X,Z;C_1 \circ
C_2)$. In each of these cases,  $C_{1} \circ C_{2}$
is again a smooth canonical relation.
However, in many interesting problems,  these assumptions  fail,
and it is important to  analyze the composition and understand
the resulting operators. It turns out that the geometry of each
canonical relation and the structure of their projections play an
important role.

\par Let $ \pi _L$ and $ \pi_ R$ be the projections, to the left
and right, from $C$ to $T^*X$ and $T^*Y$, respectively. If either
one is a local diffeomorphism, so is the other one and then $C$,
is a local canonical graph.  In the case of two canonical
relations, if at least one of $C_1$ and $C_2$ is a local canonical
graph, then $ C_1 \times C_2 $ intersects $T^*X \times
\Delta_{T^*Y}\times T^*Z$ transversally and the general
composition calculus  applies.

\par Now  consider the case when the projections are no longer local
diffeomorphisms.
When one of the projections
is singular, i.e., when the rank of its differential is
nonmaximal, then the other one is, too, and $C$ is called a {\it
singular canonical relation}. (Note: $C$ is still assumed to be smooth.)

\par Although corank($d \pi_L$)=corank($d \pi_R$) at all points, the
two  projections,
$\pi_R$ and  $\pi_L$, may have similar singularities or quite
different ones. The singularities considered in this article are
folds, submersion with folds and cross caps, which we now briefly
describe.\\

\par Let $f$ be a smooth function $f : V \to W,\ \ \dim V = \dim W=N$ and
$\mathcal{S}:={\mathcal S}(f) =\{ x \in V : \det (df (x))=0\}$.

\begin{definition}  $f$ {\em has a }  (Whitney) fold {\em
singularity along} $ {\cal S}$ {\em if} $d(det(df)) \neq 0$ {\em
on} ${\cal S}$, {\em so that}  $ {\cal S}$ {\em  is a smooth
hypersurface, } $df$ {\em drops rank by 1 there, and }
   ${\rm  Ker} \ df(x)$ {\em  intersects} $T_x{\cal S}$ {\em transversally for
every} $x \in {\cal S}$.
\end{definition}

\par Any map which has a fold singularity can be put into a local normal
form: $f(x_1,x_2, \dots, x_N)=(x_1,x_2, \dots, x_{N-1},x_N^2)$ with
respect to suitable local coordinates in the domain and codomain
\cite{gogu}.

\par Whitney folds are the singularities denoted by $S_{1,0}$ (in the
Thom theory
\cite{gogu})  and by $\Sigma_{1,0}$ (in the Boardman-Morin theory
\cite{mo1,mo2}) in the equidimensional case. The non-equidimensional
versions of
Whitney folds are submersions with folds and cross caps. We note the
difference in
notation:   when
$\dim V\leq\dim W$, the singularity classes
$S_{r,0}=\Sigma_{r,0}$, while,  if $\dim V >  \dim W$, then
$S_{r,0}=\Sigma_{r+k,0}$, where
$k=\dim V-\dim W$.

\par Let $f$ be a smooth function $f : V \to W$, dim $V = N$,  dim $W
= M$ , $ N > M.$

\begin{definition} {\em The map} $f$ {\em is a} submersion  with
folds {\em if the only singularities of} $f$  {\em are of type}
$S_{1,0}$, {\em i.e.,
of type } $\Sigma_{N-M+1,0}$.
\end{definition}

\par One checks that $f$ is a submersion with folds as follows.
At  points where
\noindent$\rank df\ge M-1$, by \cite{mo2}, we can choose suitable adapted
local coordinates on $V$ and $W$ such that $f$
has the form:
$f(x_1, x_2, \dots, x_{M-1}, x_M, \dots, x_N)=(x_1, x_2, \dots
x_{M-1}, f_1(x))$. The set $\mathcal{S}_1(f)$ where $f$ drops rank by $1$ is
described by $\mathcal{S}_1(f)= \{ x: \frac{\partial f_1}{\partial
x_i}=0, \ M \leq i
\leq N \}$. Then $f$ is a submersion with folds if  $\mathcal{S}_1(f)$ is a
smooth submanifold, i.e., $\left\{ d\left(\frac{\partial f}{\partial
x_i}\right):  \ \
M \leq  i  \leq N)\right\}$,  is linearly
independent, and if the $(N-M+1)$-dimensional kernel of
$df$ is transversal to the tangent space to $\mathcal{S}_1(f)$ in $TV$.
These conditions can be combined \cite{mo1} into
   \begin{equation}\label{swf}
\det\left[\frac{\partial^2f_1 }{\partial x_i
\partial x_j }\right]_{ M \leq i,j \leq N} \ne 0.
\end{equation}
and this is independent of the choice of adapted coordinates.

\par There are a finite number of local  normal forms for a submersion
with folds, determined by the signature of the Hessian of $f$ \cite{gogu}:

\[f(x_1,x_2, \dots, x_N) =(x_1,x_2, \dots, x_{M-1}, x_M^2 \pm
x_{M+1}^2 \pm \cdots \pm x_N^2 ).
\]
In the case relevant here, $N=M+1$ and the last
entry is a quadratic form in two variables, which is either sign
definite or indefinite; we refer to these two possibilities as
${\it elliptic}$ and ${\it hyperbolic}$ respectively.

\par We now define the second singularity class of interest; like the class
of submersions with folds, it is stable under small $C^2$ perturbations. It is
now assumed that $\dim V=N, \ \ \dim W=M$ with $ N < M$, and
$g : V\to W$ is a smooth function.

\begin{definition} {\em We say that} $g$ {\em is a} cross cap  {\em
if the only singularities of} $g$  {\em are of type} $ S_{1,0}$, {\em
i.e., of type}
$\Sigma_{1,0}$.
\end{definition}

\par To identify a cross cap, we use the description of \cite{mo1}.
At a point where  $dg$ has rank $\ge N-1$, we can find suitable
adapted coordinates such that $g(x_1, x_2, \dots, x_{N-1}, x_N)=(x_1, x_2,
\dots, x_{N-1}, g_1, g_2, \dots g_q)$, where $q= M-N+\nolinebreak1$.
The set $\mathcal{S}_1(g)$ where $g$ drops rank by $1$
is given by \linebreak$\mathcal{S}_1(g)
=\{x: \frac{\partial g_i}{\partial x_N}=0,  \quad  1 \leq i \leq q\}$.
Assume
that there is  an $i_0$, such that $\frac{\partial ^2
g_{i_0}}{\partial x_N^2} (0) \neq 0$. Then, $g$ has a cross cap
singularity near 0 if the map $\chi : \R^N \rightarrow \R^q $ given by
$\chi(x_1, x_2, \dots x_N)=(\frac{\partial g_1}{\partial x_N},
\frac{\partial g_2}{\partial x_N}, \dots, \frac{\partial
g_q}{\partial x_N})$ satisfies $\rank d\chi(0) =q$. (Notice that this forces
$N \geq q$, i.e., $M \leq 2N-1$.) These conditions can be
reformualted as: $(i)$
$\mathcal{S}_1(g)$ is smooth and of codimension $q$; $(ii)$ the
$N\times N$ minors of
$dg$ generate the ideal of $\mathcal{S}_1(g)$; and $(iii)$ $\Ker (dg)\cap
T\mathcal{S}_1(g)=(0)$.

\par  As for folds, there is a local normal form for  cross caps, due to
\linebreak Whitney\cite{wh} and Morin\cite{mo1}:

\begin{equation}\label{ccnf}
g(x_1,x_2, \dots, x_N) =(x_1,x_2, \dots, x_{N-1}, x_1 x_N, \dots
x_{M-N} x_N, x_N^2).
\end{equation}

\section{Distributions and operators associated  \\ to two cleanly intersecting
Lagrangians}

\par Classes of distributions associated to two cleanly intersecting Lagrangian
manifolds were introduced by Melrose and Uhlmann
\cite{meuh} and Guillemin and Uhlmann\cite{guuh}.  We  briefly review their
definitions and properties.

\par First, one proves that any two pairs of cleanly intersecting
Lagrangian submanifolds are (micro)locally equivalent. Thus, one can
consider the
model pair
$(\tilde{\Lambda}_0,\tilde{\Lambda}_1)$ where
$\tilde{\Lambda}_0=T_0^*{\R}^n=\{(x,\xi): x=0\}$  and
$\tilde{\Lambda}_1= N^*\{x'' = 0 \}=\{ (x,\xi): x''=\xi'=0\}$ with
$x'=(x_1, x_2,\dots, x_k)$, and $x''=(x_{k+1}, x_{k+2},\dots,
x_n)$. One defines a class of
distributions
given by oscillatory  integrals whose amplitudes  are  called {\it
product-type} symbols. Let $z=(x,s)$ be coordinates in ${\R}^m={\R}^n
\times {\R}^k$
and $(\xi, \sigma)$ the dual coordinates.

\begin{definition}$S^{p,l}(m,n,k)$ {\em  is the set of all
functions} $a(z,\xi,\sigma) \in C^{\infty} ({\R}^m \times {\R}^n
\times {\R}^k )$ {\em such that for every}
$K \subset\subset {\R}^m$ {\em and
every} $\alpha \in {\mathbb Z}^n_+, \beta \in {\mathbb Z}^k_+, \gamma \in
{\mathbb Z}^m_+$ {\em there is a }  $c_{\alpha\beta\gamma K}<\infty$
{\em such that}

\begin{equation}\label{prod-type}
|\partial_{\xi}^{\alpha}\partial_{\sigma}^{\beta}\partial_z^{\gamma}
a(z,\xi,\sigma)|
\le c_{\alpha\beta\gamma K}(1+ |\xi|)^{p- |\alpha|}  (1+ |\sigma|)^{l-|
\beta|},  \forall (z,\xi,\tau) \in K \times {\R}^n \times
{\R}^k.
\end{equation}
\end{definition}

\begin{definition} {\em \cite{guuh} Let} $I^{p,l}({\R}^n;\tilde{\Lambda}_0,
\tilde{\Lambda}_1)$ {\em be the set of all distributions } $u$
{\em such that} $u=u_1 + u_2$ {\em with} $u_1 \in C^{\infty}_0$
{\em and}
$$u_2(x)=\int e^{i((x'-s)\cdot \xi'+x'' \cdot \xi''+ s
\cdot \sigma)} a((x,s),\xi,\sigma)d\xi d\sigma ds$$
{\em with} $a \in
S^{p',l'}(m,n,k)$ {\em where} $p'=p-\frac{n}{4}+\frac{k}{2}$ {\em and}
$l'=l-\frac{k}{2}$.
\end{definition}

\par At this point, if $X$ is a manifold of dimension $n$, we can define the
class $I^{p,l}(X;\Lambda_0,
\Lambda_1)$  for any pair of Lagrangians in   $T^*X \setminus 0$ cleanly
intersecting in codimension $k$. The oscillatory integrals we use
are oscillatory integrals in sense of H\"{o}rmander[14, p.88].
\begin{definition} \emph{\cite{guuh}} $  \ u \in I^{p,l}(X;
\Lambda_0, \Lambda_1)$
{\em if} $ u=u_1 +
u_2 + \sum v_i$ {\em where} $u_1 \in I^{p+l}(\Lambda_0
   \setminus \Lambda_1)$,  $u_2 \in I^{p}(\Lambda_1 \setminus \Lambda_0)$,
{\em the sum} $\sum v_i$ {\em is locally finite and  } $v_i=Fw_i$
   {\em where} $F$ {\em  is a zero order FIO associated to} $\chi ^{-1}$
{\em where } $\chi: T^*X\setminus 0 \to T^*{\R}^n \setminus 0$
{\em  is a  canonical
   transformation such that} $\chi(\Lambda_j)\subseteq\tilde\Lambda_j, j=0,1$,
{\em microlocally, and }
$w_i
\in I^{p,l}({\R}^n;\tilde {\Lambda}_0, \tilde{\Lambda}_1)$.
\end{definition}

We say
that a distribution $u \in I^{r}(X;\Lambda_0 \setminus \Lambda_1)$
if, microlocally
away from $\Lambda_1$, $u \in I^r( X; \Lambda_0)$, the standard
H\"ormander class of
Fourier integral distributions on $X$ associated with $\Lambda_0$.

\begin{remark}\label{old-rem3.4} \emph{\cite{guuh}} $ \
{\rm If} \  u \in I^{p,l}(X;\Lambda_0, \Lambda_1)  \ {\rm then} \
   u \in  I^{p+l}(X;\Lambda_0 \setminus \Lambda_1)  \ {\text{\rm and also }}
\linebreak   u
\in I^{p}(X;\Lambda_1 \setminus \Lambda_0)$.
\end{remark}

\par We will also use the notion of   nondegenerate
phase functions which  \\ parametrize two cleanly intersecting
Lagrangians, introduced by  Mendoza [15]. Let $\lambda_0 \in
\Lambda_0 \cap \Lambda_1$ and $\Gamma \subset X \times {\R^k}
\times \left({\R}^N \setminus 0\right)$ an open, conic set.

\begin {definition} {\em \cite{men}\label{old-def3.6}  A phase function}
$\phi(x,s,\theta)$ {\em  defined
on} $\Gamma$ {\em is a } \linebreak parametrization {\em for the pair }
$(\Lambda_0,\Lambda_1)$ {\em if } \\
\indent i) $\phi_0(x,\theta):=\phi(x, 0, \theta)$ {\em where}
$\phi_0$ {\em is a nondegenerate phase function
parametrizing} $\Lambda_0$ {\em near} $\lambda_0$; {\em and} \\
\indent ii) $\phi_1(x,(\theta,\sigma)):=\phi(x,\frac{\sigma}{\mid \theta
\mid},\theta)$
{\em is  a nondegenerate phase function parametrizing}
$\Lambda_1$ {\em near} $\lambda_0$.

{\em We also refer to } $\phi_1(x;\theta;\sigma)$ {\em as a } multi-phase
function {\em for }
$(\Lambda_0,\Lambda_1)$.
\end{definition}

\par For simplicity, we now focus on the case of codimension 1
intersection relevant
here,  i.e., $k=1$. Let us consider the following example:
If $\tilde{\Lambda}_0=N^*\{ x'=0\}$ and
$\tilde{\Lambda}_1=N^*\{x=0\}$, with $x'=(x_2,x_3,\dots,x_n)$,
then $\phi(x,s,\theta')=x'\cdot \theta'\linebreak+x_1 s \mid \theta' \mid $
is a parametrization for $(\tilde{\Lambda}_0,\tilde{\Lambda}_1)$
since $\phi_0(x,\theta'):=\phi(x,0,\theta')=x'\cdot \theta'$, which is a
parametrization for $\tilde {\Lambda}_0$, and
$\phi_1(x,(\theta',\sigma)):=\phi(x,\frac{\sigma}{\mid \theta \mid},\theta)=
x'\cdot
\theta'+x_1 \sigma$ is a parametrization for $\tilde{\Lambda}_1$.

\begin{proposition}\label{old-prop3.7}  {\em \cite{men} Let} $p_1$
{\em be a homogenous function of degree} $1$ {\em such
that} $p_1(\lambda_0)=0$ {\em and} $H_{p_1}$ {\em (the Hamiltonian
vector field associated to} $p_1)$ {\em is not tangent to}
$\Lambda_0$. {\em If } $\Lambda_1$ {\em is the flow out from}
$\Lambda_0 \cap \{p_1=0\}$ {\em by } $H_{p_1}$ {\em then there is
a parametrization} $\phi$ {\em for } $(\Lambda_0,\Lambda_1)$
{\em which can be chosen such that} $\frac{\partial
\phi}{\partial s}(x,s,\theta)=p_1(x,d_x \phi)$ {\em and}
$\phi(x,0,\theta)=\phi_0$ {\em with } $\phi_0$  {\em a specified
parametrization for} $\Lambda_0$.
\end{proposition}

\begin{remark}\label{old-rem3.8} {\em [15] If}
$\phi(x;\theta; \sigma)$  \ {\em is a  multi-function for} $(\Lambda_0,
\Lambda_1)$ {\em near } $\lambda_0\in\Lambda_0\cap\Lambda_1$, {\em then, if }
$u\in\mathcal
D'(X)$ {\em has } $WF(u)$ {\em contained in a conic neighborhood of }
$\lambda_0$,
{\em then } $u\in I^{p,l}(X;\Lambda_0, \Lambda_1)$
{\em iff }
$$u(x)=\int e^{i
\phi(x;\theta;\sigma)} a(x;\theta;\sigma) d \sigma d \theta,$$
{\em  where}  $a \in S^{\tilde{p},\tilde{l}}(X\times (\R^N\setminus
0)\times\R)$,
{\em the space of } symbol-valued symbols of order
$\tilde{p},\tilde{l}$, {\em defined
by: for all }
$K\subset\subset X,\alpha\in\mathbb Z_+^N,\beta\in\mathbb Z_+,
\gamma\in\mathbb Z_+^n,$

\begin{equation}\label{ests-svs}|\partial^\alpha_{\theta}
\partial^\beta_{\sigma}\partial^\gamma_{x}
a(x;\theta;\sigma)| \leq
c_{\alpha\beta\gamma
K}(1+|\theta|+|\sigma|)^{\tilde{p}-|\alpha|}(1+|\sigma|)^{\tilde{l}-\beta},
\end{equation}
{\em with } $p=\tilde{p}+\tilde{l}+\frac{N+1}{2}-\frac{n}{4},
l=-\tilde{l}-\frac{1}{2}$.

\end{remark}

Finally, we define the classes of generalized (or {\it paired
Lagrangian}) Fourier
integral operators, to one of which we will show the normal operator
$F^*F$ belongs.
Recall that a {\em canonical relation} $C\subset(T^*X\setminus 0)\times
(T^*Y\setminus 0)$ is a smooth, conic submanifold such that
$C':=\left\{(x,y;\xi,\eta): (x,\xi;y,-\eta)\in C\right\}$ is a
Lagrangian submanifold
of $\left(T^*(X\times Y),\omega_{T^*(X\times Y)}\right)$, i.e,  $C$
is a Lagrangian
with respect to $\omega_{T^*X}-\omega_{T^*Y}$.

\begin{definition}\label{old-def3.5} {\em If } $C_0,C_1\subset
(T^*X\setminus 0)\times
(T^*Y\setminus 0)$ {\em are smooth, conic canonical relations
intersecting cleanly,
then } $I^{p,l}(C_0,C_1):=I^{p,l}(X,Y;C_0,C_1)$  {\em denotes the set 
operators }
$F:{\cal E}'(Y)\longrightarrow{ \cal D}'(X)$ {\em whose Schwartz
kernels are in }
$I^{p,l}(X\times Y; C_0',C_1')$.
\end{definition}

\section{Fold caustics in the marine geometry}

\par In three spatial dimensions, let $s$ be a fixed source on the
surface,\linebreak
$\{x_3=0\}$;
$H(x,\xi)=\frac{1}{2}(c_0(x)^{-2}-|\xi|^2)$  the Hamiltonian
associated to the smooth
background soundspeed $c_0(x)$ in (\ref{starstar}); and $\Lambda_s$ 
the image of
$T_s^*\R^3
\setminus 0$ under the bicharacteristic flow  associated to $H$,
which is a Lagrangian submanifold of $T^*\R^3 \setminus 0$. The
assumption of a (point) fold caustic means that the only singularities of
the spatial projection $\pi_Y : \Lambda_s \rightarrow Y$ are folds. We make
use of the description of $\Lambda_s$ in Nolan\cite{no}. It can be
parametrized by $t_{inc}$, the time travelled by the incident
ray, and the takeoff direction $(p_1,p_2,p_3) \in \mathbb S ^2$. We
can change these coordinates to $(x_1, x_2, p_3)$ \cite{no}. Hence on
$\Lambda_s$, $x_3=f(x_1, x_2, p_3)$ and $(p_1,p_2)=(g_1(x_1, x_2,
p_3), g_2(x_1, x_2, p_3))$. In this new setting,  det $d\pi_S= \frac
{\partial f}{\partial p_3}(x_1, x_2, p_3)$ and  fold
caustics occur where $ \frac{\partial f}{\partial p_3} = 0$
and $\frac{\partial^2 f}{\partial p_3^2} \neq 0$.

\par In the marine geometry, the source $s=(s_1, s_2,0)$ is subject
to the restriction $s_2=r_2$, so we
consider  just $s_1$ as an independent coordinate. Fix an $s_2\in\R$,
i.e., consider a single pass of the vessel, and let $\Lambda_{s_2}$
be the union of
the flowouts
$\left\{
\Lambda_{(s_1,s_2,0)}: s_1\in\R\right\}$, so that $\Lambda_{s_2}\subset
T^*\R^3\setminus 0$ is an involutive submanifold.
We say that fold
caustics (and no
worse) appear for the background sound speed if, considering
$s_1$ as a variable, the spatial projection
$\pi_Y:\Lambda_{s_2}\longrightarrow \R^3$ is a submersion with folds. By the
structural stability of submersions with folds, this condition will
then hold for all
$s_2'$ close to $s_2$. The presence of fold caustics may be
characterized as follows.
The variables
$x_3$ and
$(p_1, p_2)$ are functions of the other four:
$x_3=f(x_1, x_2, s_1, p_3)$,  $(p_1,p_2)=(g_1(x_1, x_2,s_1, p_3),
g_2(x_1, x_2, s_1, p_3))$. The differential $d \pi_Y$ then becomes:

\[{d \pi_Y}= \left(\begin{array}{cccc}
1 & 0 & 0 & 0 \\
0 & 1 & 0 & 0 \\
\frac{\partial f} {\partial x_1} & \frac{\partial f}{\partial x_2}
& \frac{\partial f}{\partial s_1} & \frac{\partial f}{\partial
p_3}
\end{array} \right).\]

We have

\begin{eqnarray*}
{\rm rank}  \ \  d \pi_Y = \left\{
\begin{array}{ccc}
2, & {\rm if}&  \frac{\partial f}{\partial s_1}=\frac{\partial
f}{\partial p_3}=0\\
3,  & {\rm if }&   \frac{\partial f}{\partial s_1} \neq 0   \
{\rm or} \ \frac{\partial f}{\partial p_3} \neq 0.
\end{array}  \right.
\end{eqnarray*}

\par Suppressing $s_2$, let $ {\mathcal S}_1^{\Lambda}:=
{\mathcal S}_1(\pi_Y)=\{  \frac{\partial f}{\partial p_3}=
\frac{\partial f}{\partial s_1}=0 \}$ be the critical set of $\pi_Y$, where
$\textrm{ rank } d\pi_Y$
drops by 1. At points of
${\mathcal S}_1^{\Lambda}$,
$\Ker d\pi_Y$ is two-dimensional and spanned by
$\{ (0,0, \delta s_1, \delta p_3)  \}$.  The tangent space to ${\cal
S}_1^{\Lambda}$ is
\begin{equation}\label{tansone}
T {\cal S}_1^\Lambda= \Ker (d_{x_1,x_2,p_3,s_1}(\frac{\partial
f}{\partial p_3})) \cap \Ker (d_{x_1,x_2,p_3,s_1}(\frac{\partial
f}{\partial s_1})).
\end{equation}
where

\begin{equation}\label{dpthree}
d_{x_1,x_2,p_3,s_1}(\frac{\partial
f}{\partial p_3})= (\frac{\partial^2 f}{\partial x_1 \partial
p_3}, \frac{\partial^2 f}{\partial x_2
\partial p_3}, \frac{\partial^2 f}{
\partial p_3^2}, \frac{\partial^2 f}{\partial s_1 \partial p_3})
\end{equation}
and

\begin{equation}\label{dsone}
d_{x_1,x_2,p_3,s_1}(\frac{\partial
f}{\partial s_1})= (\frac{\partial^2 f}{\partial x_1 \partial
s_1}, \frac{\partial^2 f}{\partial x_2
\partial s_1}, \frac{\partial^2 f}{
\partial p_3 \partial s_1}, \frac{\partial^2 f}{\partial s_1 ^2}).
\end{equation}
Then, $\pi_Y$ is a submersion with folds if

\begin{eqnarray}\label{starstarstar}
{\cal S}_1^{\Lambda} \textrm{ is smooth, i.e., the
gradients in (\ref{dpthree}) and (\ref{dsone}) are}\ \ \ \ \nonumber\\
\textrm{linearly independent, and }T{\mathcal S}_1^{\Lambda}\textrm{
is transversal
to }\Ker d\pi_S,\
\end{eqnarray}
i.e., if
\begin{equation}\label{detnz}
\left|\begin{array}{cc}
\frac{\partial^2 f}{\partial p_3^2}& \frac{\partial^2 f}{\partial s_1
\partial p_3}\\
\frac{\partial^2 f}{\partial p_3 \partial s_1}& \frac{\partial^2
f}{\partial s_1 ^2}
\end{array}\right|\ne 0.
\end{equation}

\par Next, we  parametrize the canonical relation $C$ of $F$ in
terms of $s_1, x_1, x_2$ and
$p_3$ ; $(\alpha_1, \alpha_2, \sqrt {1-|\alpha|^2})$, the take off
direction of the reflected ray, writing $\alpha=
(\alpha_1,\alpha_2)$; and $\tau$, the variable dual  to time.
Following \cite{no}, the canonical relation $C \subset T^*(\Sigma_{r,s} \times
(0,T)) \times T^*\R^3_+$ is parametrized as

\begin{eqnarray*}
C = \big\{ (&s_1&, r_1(\cdot), r_2(\cdot),
t_{inc}(\cdot) + t_{ref}(\cdot),
      \sigma(\cdot),
\rho_1(\cdot),   \rho_2(\cdot), \tau ;  \\
   &{ } x_1&, x_2, f(\cdot), -\tau(c_0^{-1}(\cdot)
\alpha_1 + g_1(\cdot)),  -\tau (c_0^{-1}(\cdot)\alpha_2 + g_2(\cdot)),\\
    &{ }&-\tau(c_0^{-1}(\cdot), \alpha)\sqrt{1-|\alpha|^2} + p_3) \big\}
\end{eqnarray*}
where
\begin{eqnarray*}
f(\cdot)&=&f(x_1,x_2, s_1, p_3); \\
r_j(\cdot)&=&r_j(x_1,x_2, f(x_1,x_2,
s_1, p_3),  j=1,2;\\
   t_{inc}(\cdot)&=&t_{inc}(x_1, x_2, p_3); \\
t_{ref}(\cdot)&=&t_{ref}(x_1,x_2, f(x_1,x_2, s_1, p_3);\\
\sigma(\cdot)&=&\sigma(x_1,x_2, f(x_1,x_2, s_1, p_3);\\
\rho_j(\cdot)&=&\rho_j(x_1,x_2, f(x_1,x_2, s_1, p_3), \alpha),  j=1,2;\\
g_j(\cdot)&=&g_j(x_1, x_2, s_1, p_3),  j=1,2 {\rm ; and }\\
c_0^{-1}(\cdot)&=&c_0^{-1}\left(x_1, x_2,f(\cdot)\right).
\end{eqnarray*}

It was proven in \cite{nosy} that $F:\mathcal E'\left(\mathbb
R_+^3\right)\longrightarrow\mathcal D'\left(\Sigma_{r.s}\times
(0,T)\right)$ is a
Fourier integral operator,
$F\in I^{\frac{3}{4}} (\Sigma_{r,s} \times (0,T), \R^3_+; C)$; we now
show that the presence of  caustics of fold type (and no worse) imposes certain
conditions on $C$, namely that $\pi_R : C \rightarrow T^*\R^3
\setminus 0$ is a submersion with folds and  \linebreak $\pi_L : C
\rightarrow T^*(\Sigma_{r,s} \times (0,T)) \setminus 0$ is a cross
cap. In fact, with respect with the coordinates above,
$\pi_R :
\R^7 \rightarrow \R^6$ is given by

\begin{eqnarray*}
\pi_R(x_1,x_2,p_3,s_1,\alpha_1,\alpha_2,\tau)= \big ( x_1, x_2,
f(\cdot); & -&\tau(c_0^{-1}(\cdot) \alpha_1 + g_1(\cdot)),\\
&-&\tau(c_0^{-1}(\cdot)
\alpha_2 + g_2(\cdot)),\\
&-&\tau(c_0^{-1}(\cdot) \sqrt{1-|\alpha|^2} + p_3)
\big ).
\end{eqnarray*}
Thus,

\[{d \pi_R}= \left(\begin{array}{ccccccc}
1 & 0 & 0 & 0 & 0 & 0 & 0\\
0 & 1 & 0 & 0 & 0 & 0 & 0\\
\frac{\partial f}{\partial x_1} & \frac{\partial f}
{\partial x_2} & \frac{\partial f}{\partial p_3} & \frac{\partial
f}{\partial s_1} & 0 & 0 & 0\\
A_1 & A_2 & A_3 & A_4 & -\tau c_0^{-1} & 0 & -(c_0^{-1}\alpha_1+g_1)\\
B_1 & B_2 & B_3 & B_4 & 0 & -\tau c_0^{-1} & -(c_0^{-1}\alpha_2+g_2)\\
C_1 & C_2 & C_3 & C_4 &
\frac{-\tau c_0^{-1}\alpha_1}{\sqrt{1-|\alpha|^2}} &
\frac{-\tau c_0^{-1}\alpha_2}{\sqrt{1-|\alpha|^2} } &
   -(c_0^{-1}\sqrt{1-|\alpha|^2}+p_3)
\end{array} \right)\]

\par Thus,  rank $d \pi_R = \left\{ \begin{array}{ccc}
5, & {\rm if} &  \ \ \frac{\partial f}{\partial p_3}= \frac{\partial
f}{\partial s_1}=0\\
6,  & {\rm if } &  \ \ \frac{\partial f}{\partial p_3} \neq 0   \
{\rm or} \ \frac{\partial f}{\partial s_1} \neq 0
\end{array} \right. $ because the matrix

\[ \left(\begin{array}{ccc}

   -\tau c_0^{-1} & 0 & -(c_0^{-1}\alpha_1+g_1)\\

   0 & -\tau c_0^{-1} & -(c_0^{-1}\alpha_2+g_2)\\

   - c_0^{-1} \frac{\tau\alpha_1}{\sqrt{1-|\alpha|^2}} & -c_0^{-1}

\frac{\tau\alpha_2}{\sqrt{1-|\alpha|^2} } & -(c_0^{-1}\sqrt{1-|\alpha|^2}+p_3)

\end{array} \right)\]
is nonsingular \cite{no}. Hence, the critical set ${\mathcal S}_1^C:= {\mathcal
S}_1(\pi_R)$
   is a smooth, codimension two submanifold.
   (Recall that by general considerations
[4], this must equal ${\mathcal S}_1(\pi_L)$, and $d\pi_R$ and
$d\pi_L$ must drop rank
by the same amount at each point.) At these points,
$\Ker d
\pi_R =
\{ (0, 0,
\delta p_3,
\delta s_1,
\delta \alpha_1,
\delta \alpha_2, \delta \tau) \}$ where $\delta \alpha_1, \delta
\alpha_2, \delta \tau$ depend on $\delta p_3, \delta s_1$.
The tangent space to ${\cal S}_1^C$ is

\[
T {\cal S}_1^C= \Ker
\left(d_{x_1,x_2,p_3,s_1,\alpha_1,\alpha_2,\tau}(\frac{\partial
f}{\partial p_3})\right) \cap \Ker
\left(d_{x_1,x_2,p_3,s_1,\alpha_1,\alpha_2,\tau}(\frac{\partial
f}{\partial s_1})\right),
\]
with

\begin{equation}\label{grad-one}
d_{x_1,x_2,p_3,s_1,\alpha_1,\alpha_2,\tau}(\frac{\partial
f}{\partial p_3})= (\frac{\partial^2 f}{\partial x_1 \partial
p_3}, \frac{\partial^2 f}{\partial x_2
\partial p_3}, \frac{\partial^2 f}{
\partial p_3^2}, \frac{\partial^2 f}{\partial s_1 \partial p_3}, 0, 0,0)
\end{equation}
and

\begin{equation}\label{grad-two}
d_{x_1,x_2,p_3,s_1,\alpha_1,\alpha_2,\tau}(\frac{\partial
f}{\partial s_1})= (\frac{\partial^2 f}{\partial x_1 \partial
s_1}, \frac{\partial^2 f}{\partial x_2
\partial s_1}, \frac{\partial^2 f}{
\partial p_3 \partial s_1}, \frac{\partial^2 f}{\partial s_1 ^2}, 0, 0,0),
\end{equation}
and so we  see that  $\Ker d \pi_R$ is transversal to $T {\cal S}_1^C$
because of the  condition (\ref{starstarstar}). Hence, $\pi_R$ is a submersion
with folds. Note that without further restrictions on $f$, the
projection $\pi_R$
can either an elliptic or hyperbolic submersion with folds.

\par Similarly, with respect to the above coordinates,
reordered for ease of display, $\pi_L: \R^7
\rightarrow \R^8$ is given by

\[\pi_L(s_1, x_1, x_2, \alpha_1,\alpha_2, p_3, \tau)= \big(s_1,
r_1(\cdot),   r_2(\cdot),
t_{inc}(\cdot) + t_{ref}(\cdot); \sigma(\cdot),
\rho_1(\cdot), \rho_2(\cdot), \tau \big),
\]
and thus ${d \pi_L}= $

\[
\left(\begin{array}{ccccccc}
1 & 0 & 0 & 0 & 0 & 0 & 0\\
\frac{\partial r_1}{\partial s_1} & \frac{\partial r_1}{\partial
x_1} +  \frac{\partial r_1}{\partial x_3}\frac{\partial
f}{\partial x_1} & \frac{\partial r_1}{\partial x_2} +
\frac{\partial r_1}{\partial x_3}\frac{\partial f}{\partial x_2} &
\frac{\partial r_1}{\partial \alpha_1} & \frac{\partial r_1}{\partial
\alpha_2} & \frac{\partial r_1}{\partial x_3} \frac{\partial f}{\partial
p_3} & 0\\
\frac{\partial r_2}{\partial s_1} & \frac{\partial r_2}{\partial
x_1} + \frac{\partial r_2}{\partial x_3}\frac{\partial f}{\partial
x_1} & \frac{\partial r_2}{\partial x_2}  +  \frac{\partial
r_2}{\partial x_3}\frac{\partial f}{\partial x_2}& \frac{\partial
r_2}{\partial \alpha_1}
& \frac{\partial r_2}{\partial \alpha_2} & \frac{\partial r_2}{\partial
x_3} \frac{\partial f}{\partial p_3} & 0\\
\frac{\partial t_{ref}}{\partial s_1} & \frac{\partial
t_{inc}}{\partial x_1} + \frac{\partial t_{ref}}{\partial
x_3}\frac{\partial f}{\partial x_1} & \frac{\partial
t_{inc}}{\partial x_2} +  \frac{\partial t_{ref}}{\partial
x_3}\frac{\partial f}{\partial x_1} & \frac{\partial
t_{ref}}{\partial \alpha_1} &
\frac{\partial t_{ref}}{\partial \alpha_2} & \frac{\partial
t_{inc}}{\partial p_3} + \frac{\partial t_{ref}}{\partial x_3}
\frac{\partial f}{\partial p_3} & 0 \\
\frac{\partial \sigma}{\partial s_1} & \frac{\partial
\sigma}{\partial x_1} +  \frac{\partial \sigma}{\partial
x_3}\frac{\partial f}{\partial x_1} & \frac{\partial
\sigma}{\partial x_2} +  \frac{\partial \sigma}{\partial
x_3}\frac{\partial f}{\partial x_2}& \frac{\partial
\sigma}{\partial \alpha_1} & \frac{\partial \sigma}{\partial
\alpha_2} & \frac{\partial \sigma}{\partial x_3}
\frac{\partial f}{\partial p_3} & 0 \\
\frac{\partial \rho_1}{\partial s_1} & \frac{\partial
\rho_1}{\partial x_1} +  \frac{\partial \rho_1}{\partial
x_3}\frac{\partial f}{\partial x_1}& \frac{\partial
\rho_1}{\partial x_2}  +  \frac{\partial \rho_1}{\partial
x_3}\frac{\partial f}{\partial x_2} & \frac{\partial
\rho_1}{\partial \alpha_1} & \frac{\partial \rho_1}{\partial
\alpha_2} &
\frac{\partial \rho_1}{\partial x_3} \frac{\partial \rho_1}{\partial p_3}
& 0\\
\frac{\partial \rho_2}{\partial s_1} & \frac{\partial
\rho_2}{\partial x_1} +  \frac{\partial \rho_2}{\partial
x_3}\frac{\partial f}{\partial x_1}& \frac{\partial
\rho_2}{\partial x_2} +  \frac{\partial \rho_2}{\partial
x_3}\frac{\partial f}{\partial x_2} & \frac{\partial
\rho_2}{\partial \alpha_1} & \frac{\partial \rho_2}{\partial
\alpha_2} &
\frac{\partial \rho_2}{\partial x_3} \frac{\partial \rho_2}{\partial p_3}
& 0\\
0 & 0 & 0 & 0 & 0 & 0 & 1
\end{array}  \right)\]

\par Since $t_{inc}$ and $p_3$ are  independent coordinates, we have
$\frac{\partial t_{inc}}{\partial
p_3}=0$. Also, because of the choice of the coordinates $x_1, x_2,
x_3$, it follows that $\frac{\partial f}{\partial x_1}= \frac{\partial
f}{\partial x_2}=0$ at the caustic points \cite{no-per}.

\par From this, it follows that the rank of  $d \pi_L = \left\{
\begin{array}{ccc}
6, & {\rm if} &  \ \ \frac{\partial f}{\partial p_3}= \frac{\partial
f}{\partial s_1}=0\\
7,  & {\rm if } &  \ \ \frac{\partial f}{\partial p_3} \neq 0   \
{\rm or} \ \frac{\partial f}{\partial s_1} \neq 0,
\end{array}  \right. $  \\ because the matrix

\[ \left(\begin{array}{cccc}
\frac{\partial r_1}{\partial x_1}&   \frac{\partial r_1}{\partial x_2}
&  \frac{\partial r_1}{\partial \alpha_1} &  \frac{\partial r_1}{\partial
\alpha_2} \\
\frac{\partial r_2}{\partial x_1}&   \frac{\partial r_2}{\partial x_2}
&  \frac{\partial r_2}{\partial \alpha_1} &  \frac{\partial r_2}{\partial
\alpha_2} \\
\frac{\partial \rho_1}{\partial x_1}&   \frac{\partial \rho_1}{\partial
x_2}  &  \frac{\partial \rho_1}{\partial \alpha_1} &  \frac{\partial
\rho_1}{\partial
\alpha_2}\\
\frac{\partial \rho_2}{\partial x_1}& \frac{\partial
\rho_2}{\partial x_2}  &  \frac{\partial \rho_2}{\partial
\alpha_1} & \frac{\partial \rho_2}{\partial \alpha_2}
\end{array} \right)\]
is nonsingular \cite{no}. Furthermore, $\mathcal{S}_1(\pi_L)=\mathcal{S}_1^C=
\{\frac{\partial f}{\partial p_3}=\frac{\partial f}{\partial
s_1}=0\}$ is smooth and
the $7\times7$ minors of $d\pi_L$ generate the ideal of $\mathcal{S}_1^C$.
Finally, it also
follows that

\[\Ker d \pi_L = \{ (0, \delta x_1, \delta x_2, \delta \alpha_1,
\delta \alpha_2, \delta p_3, 0) \},
\]
where   $\delta \alpha_1,
\delta \alpha_2, \delta x_1, \delta x_2$ depend on
$\delta p_3$, and thus  $\Ker d\pi_L\cap T {\cal S}_1^C=(0)$. Hence, $\pi_L$
is a cross cap.

\ \

\par This leads us to formulate a general class of canonical
relations with this structure.
\begin {definition}\label{def-fcc}
{\rm Let}  $X$ { \rm and} $Y$ {\rm be manifolds of dimensions} $n$
{\rm and} $n-1,$ {\rm respectively, and let} $C$ {\rm be a
canonical relation in } $\left(T^*X \setminus 0\right) \times
\left(T^*Y \setminus 0\right)$. $C $  {\rm is a} {\it folded cross
cap} {\rm canonical
relation if:}

a) $\pi_R : C \rightarrow T^*Y \setminus 0$ {\rm is a submersion with
folds, with  singular set} ${\cal S}_1,$
{\rm and  } $\pi_R({\cal S}_1)$ {\rm is a nonradial hypersurface in}
$T^*Y \setminus 0$; {\rm and}\\
\indent b) $\pi_L : C \rightarrow T^*X \setminus 0$ {\rm has a cross cap
singularity along} ${\cal S}_1$ {\rm and} $\pi_L({\cal S}_1)$ {\rm
is a nonradial submanifold.}

{\rm We say that} $C$ {\rm is an
\emph{elliptic}, respectively,
\emph{hyperbolic}, folded cross cap if } $\pi_R$ {\rm is an elliptic,
respectively, hyperbolic,
submersion with folds}.
\end{definition}
\begin{remark} {\rm We will see in \S6 (see discussion following Prop.
6.2) that } $\pi_L(\so)$  {\rm is necessarily } maximally
noninvolutive{\rm , defined
after (\ref{iii}) below.}
\end{remark}

\par Here, as usual, a conic submanifold $\Gamma \subset
T^*Y\setminus 0$ is {\it
nonradial} if for all $(y,\eta)\in\Gamma$, the canonical one-form
$\sum \eta_j dy_j$
does not vanish identically on $T_{(y,\eta)}\Gamma$. Since $\pi_R(\so)$ is the
immersed image of $\so$, $\pi_R(\so)$ is nonradial at $\pi_R(c)$ iff $\sum
(\pi_R^*\eta_j) d\pi_R^*(d y_j)\ne 0$ as an element of $T_{c}^*C$. We can
understand the significance of this for the canonical relation
arising in the marine
geometry as follows. Using the fact that $df=0$ at the caustics, the
expressions for
$\pi_R$ and $d\pi_R$ computed above imply that, at $c\in\so$,
\begin{equation}\label{pb}
\sum_{j=1}^3(\pi_R^*\eta_j) d\pi_R^*(d
y_j)=-\tau\left((c_0^{-1}\alpha_1+g_1) dx_1+
(c_0^{-1}\alpha_2+g_2) dx_2\right).
\end{equation}
In order  to be 0 on $T_{c}\so$, this must be a linear combination of
(\ref{grad-one}) and (\ref{grad-two}). By (\ref{detnz}), the only possible
linear combination is the trivial one; however, by the expression for
$\pi_R$, this
forces $\eta=(0,0,\eta_3)$ for some $\eta_3\ne 0$. That is, the fold
caustic surface
in $Y$ must be horizontal at this point. While there certainly exist background
soundspeeds $c_0(\cdot)$ for which this happens, the most basic
examples of fold
caustics arising from refraction about a low velocity lens
   (see Nolan and Symes \cite{nosy-fold}) have fold
surfaces which are not horizontal.

\par Similarly, since $\pi_L|_{\so}$ is an immersion,
$\pi_L(\so)\subset T^*X$ is
nonradial iff
\linebreak$\sum (\pi_L^*\xi_j) d\pi_L^*(d x_j)\ne 0$ in $T^*_cC$. The
physical interpretation
of a radial point, i.e., a point $c\in\so$ where this fails,  for the
marine geometry
is less clear and, to proceed, we will simply need to assume that
such points are
absent.
\vfil\eject

\section{Model case}

\par We showed in the previous section that the canonical relation
$C$ arising from the marine geometry in the presence of fold caustics is a
folded cross cap.  To help understand the nature of the normal operator, we
first consider
a model folded cross cap canonical relation,
$C_0$, in
$(T^*\R^n \setminus 0) \times (T^*\R^{n-1} \setminus 0)$,
parametrized by the phase function
\[\phi_0(x,y,\theta'')=(x''-y'') \cdot \theta'' + \left((x_n^2-x_{n-1}^2)
y_{n-1}+x_ny_{n-1}^2\right)\theta_1,\]

\noindent where $(x, y, \theta'')
\in \R^n \times \R^{n-1} \times (\R^{n-2} \setminus 0)$, in the
region $\{ |\theta_1| \geq c|\theta''| \}$. Here and at various
points below we use
the notation $x=(x_1,\dots,x_n)=(x'', x_{n-1}, x_n)=(x_1,x''',x_n)$, $y=(y'',
y_{n-1})$ and
$\theta''= (\theta_1,\theta''')$.

\par For simplicity, in this section and the next one  we will make
the  choice
that $C$ is a hyperbolic folded cross cap. This corresponds to the
choice of the $(-)$ sign in the
$(x_{n-1}^2 - x_n^2)$ term of $\phi_0$. There are no significant changes
needed in the calculations  for the  elliptic case.

\par One easily calculates that

\begin{eqnarray*}
C_0=\big\{( x'', &x_{n-1}, x_{n}, \theta'', -2x_{n-1}y_{n-1}\theta_1,
y_{n-1}^2\theta_1+2x_ny_{n-1}\theta_1;\\
&y'', y_{n-1}, \theta'', -\left(2x_n
y_{n-1}+(x_n^2-x_{n-1}^2)\right)\theta_1)\\
&:
|\theta_1| \geq c|\theta''|, \  x_i=y_i,  \   2 \leq i \leq n-2,\\
&x_1-y_1
+(x_n^2-x_{n-1}^2)y_{n-1} +x_n y_{n-1}^2=0  \big\}
\end{eqnarray*}

\par We will verify that the projections of $C_0$ to the left and to
the right have the desired singularities. Note that $(x, y_{n-1},
\theta'')$ are coordinates on $C_0$. Hence, reordering the variables 
for ease of display,

\begin{eqnarray*}
\pi_R (x_1,x''', y_{n-1},\theta_1, \theta''',x_{n-1},x_n)= \big( 
x_1& +&(x_n^2-x_{n-1}^2)y_{n-1} +x_n
y_{n-1}^2, x''', y_{n-1}\\
&;& \theta'', -\left(2x_n
y_{n-1}+(x_n^2-x_{n-1}^2)\right)\theta_1\big)
\end{eqnarray*}

and

\[{d \pi_R}= \left(\begin{array}{ccccccc}
1 & 0 & A & 0 & 0 & B & D\\
0 & I_{n-3} & 0 & 0 & 0 & 0 & 0\\
0 & 0 & 1 & 0 & 0 & 0 & 0\\
0 & 0 & 0 & 1 & 0 & 0 & 0\\
0 & 0 & 0 & 0 & I_{n-3} & 0 & 0\\
0 & 0 & E & F & 0 & 2x_{n-1}\theta_1 &
-2y_{n-1}\theta_1-2x_n\theta_1
\end{array} \right)\]
where $A= (x_n^2-x_{n-1}^2)+2x_ny_{n-1}, B=-2x_{n-1} y_{n-1},
D= 2x_ny_{n-1}+y_{n-1}^2, \\ E= -2x_n\theta_1$, and  $F=
-2x_ny_{n-1}-(x_n^2-x_{n-1}^2)$. This is a $ (2n-2) \times (2n-1)$
matrix, with
$$ {\rm rank}  \ \  d \pi_R = \left\{
\begin{array}{ccc}
2n-3, & {\rm if} &  \ \ x_{n-1}=x_n+y_{n-1}=0\\
2n-2,  & {\rm if } &  \ \ x_{n-1} \neq 0   \ {\rm or} \  x_n+y_{n-1}
\neq 0.
\end{array}  \right. $$

\par Let ${\cal S}_1^{C_0}:={\cal S}_1(\pi_R)=\{ (x,y_{n-1},\theta'') :
x_{n-1}=x_n+y_{n-1}=0
\} $  be the set where $ d\pi_R$ drops rank by $1$. Off of this
smooth, codimension
two submanifold, $\pi_R$ is a submersion. The kernel of $ \pi_R$ at
${\cal S}_1^{C_0}$ is spanned by $\{ \frac{\partial {\quad
}}{\partial x_{n-1}},
\frac{\partial { \quad}}{\partial x_n} \}$ and thus  intersects the tangent
space $T {\cal S}_1^{C_0}$ transversally. We  also note that  the
Hessian, $(x_{n-1}, x_n) \rightarrow
-(x_n^2-x_{n-1}^2)\theta_1$, is sign-indefinite. We
thus conclude that $\pi_R$ is a hyperbolic submersion with folds.

\par As for $\pi_L$, we have
$$\pi_L(x'',x_{n-1},x_n, \theta_1, \theta''',y_{n-1})=\left(x;
\theta'', -2x_{n-1}y_{n-1}\theta_1,
\left(y_{n-1}^2+2x_ny_{n-1}\right)\theta_1\right),$$
so that

\[{d \pi_L}=
\left(\begin{array}{cccccc} I_{n-2} & 0 & 0 & 0 & 0 &
0\\
0 & 1 & 0 & 0 & 0 & 0 \\
0 & 0 & 1 & 0 & 0 & 0 \\
0 & 0 & 0 & 1 & 0 & 0 \\
0 & 0 & 0 & 0 & I_{n-3} & 0 \\
0 & A' & 0 & B' & 0 & -2x_{n-1}\theta_1 \\
0 & 0 & D' & E' &  0 & 2(x_n+y_{n-1})\theta_1
\end{array} \right)\]
where $A'= -2y_{n-1}\theta_1,  B'=-2x_{n-1} y_{n-1},
D'= 2y_{n-1}\theta_1$, and $E'= 2x_ny_{n-1}+\nolinebreak y_{n-1}^2$.
This is a $ 2n \times (2n-1)$ matrix and
$$ {\rm rank } \ d \pi_L =
\left\{ \begin{array}{ccc}
2n-2, & {\rm if} &  \ \ x_{n-1}=x_n+y_{n-1}=0\\
2n-1,  & {\rm if } &  \ \ x_{n-1} \neq 0   \ {\rm or} \  x_n+y_{n-1}
\neq 0.
\end{array}  \right. $$

\par Thus, $d\pi_L$ drops rank by $1$ at ${\cal S}_1^{C_0}$, and  off
of this set is an immersion. The kernel of $ \pi_L$ is spanned by $\{
\frac{\partial\quad}{\partial y_{n-1} }\}$ and intersects the tangent
space $T {\cal S}_1^{C_0}$ transversally. Also, the rank of the differential of
$(x,y_{n-1},\theta'') \rightarrow\linebreak
(-2x_{n-1}\theta_1, \left(2y_{n-1} +2 x_n\right)\theta_1)$ is $2$ at
${\cal S}_1^{C_0}$. We conclude that $\pi_L$ has a cross cap
singularity.

\par We note for future use the images
of $ {\cal S}_1^{C_0}$ under $\pi_L$ and $\pi_R$: since $\pi_R({\cal
S}_1^{C_0})=\{
(x'', -x_n, \theta'', x_n^2 \theta_1) \}$ and $ \pi_L({\cal
S}_1^{C_0})=\{ (x'', 0, x_n,\theta'', 0, -x_n^2 \theta_1) \}$, we have

\begin{equation}\label{ii}
\pi_R\left({\cal S}_1^{C_0}\right)= \left\{
\xi_{n-1}= x_n^2 \xi_1  \right\}
\end{equation}
and

\begin{equation}\label{iii}
\pi_L\left({\cal S}_1^{C_0}\right)= \left\{
x_{n-1}=\xi_{n-1}=
\xi_n+x_n^2\xi_1=0\right\}.
\end{equation}

Note that $\pi_R\left({\cal S}_1^{C_0}\right)$ is a nonradial hypersurface in
$T^*\mathbb R^{n-1}\setminus 0$, while
$\pi_L({\cal S}_1^{C_0})$ is a nonradial, codimension three submanifold
of $T^*\mathbb R^{n}\setminus 0$, given by
defining functions $p_1=\xi_{n-1}, p_2=\xi_1 x_{n-1},
p_3=\xi_n+x_n^2\xi_1$ with  Poisson brackets $\{p_1, p_2 \}=1, \{p_1,p_3 \}=0,
\{p_2,p_3 \}=0$.  This is \emph{maximally noninvolutive} in the sense that
$\omega_{T^*\R^n} |_{\pi_L({\cal S}_1)}$ has the maximal possible rank for a
codimension three submanifold of $T^*\mathbb R^{n}\setminus 0$, namely
$2n-4$.

\ \

\par Next, we  calculate the composition $C_0^t \circ C_0=\{ (x,
\xi; y, \eta): \exists (z, \zeta) \ \ {\rm such
   \ \ that}\\ \ \ (x, \xi; z, \zeta) \in C_0^t  \  {\rm and }  \
   (z, \zeta; y, \eta) \in C_0 \}$. We have $(x, \xi; z, \zeta) \in C_0^t  $
iff $(z,\zeta;x,\xi)\in\nolinebreak C_0$ iff

$$  z_i=x_i, \ \  2 \leq i
   \leq n-2;$$

$$ z_1-x_1 +(z_n^2-z_{n-1}^2)y_{n-1} +z_n
x_{n-1}^2=0; $$

$$   \zeta_i=\xi_i   \ \  1 \leq i  \leq n-2; $$

$$ \zeta_n =2 z_n x_{n-1} \xi_1 + x_{n-1}^2\xi_{1}; $$

$$ \zeta_{n-1}= -2z_{n-1} x_{n-1} \xi_1 ;\textrm{ and } $$

$$ \xi_{n-1} = -(z_n^2-z_{n-1}^2)\xi_1-2z_n x_{n-1} \xi_1,$$
with $(z, \zeta; y, \eta) \in C_0 $ being determined by the same
equations, but where
$(y,\eta)$ replaces $(x,\xi)$. We thus obtain that $C_0^t\circ C_0$ consists of
$(x,\xi;y,\eta)$ such that, for some $(z_{n-1},z_n)\in\mathbb R^2$,

$$x_i=y_i, \ \ 2 \leq i \leq n-2; $$

$$ \xi_i=\eta_i,  \ \  1 \leq i \leq n-2 ; $$

$$ z_{n-1}(x_{n_1}-y_{n-1})\xi_1=0 ;$$

$$ (x_{n-1}-y_{n-1})(2 z_n +x_{n-1}+ y_{n-1}) =0 ;$$

$$ \xi_{n-1} = -\left((z_n^2-z_{n-1}^2)+2z_n x_{n-1}\right)
\xi_1;$$

$$ \eta_{n-1} = -\left((z_n^2-z_{n-1}^2)+2z_n y_{n-1}\right) \eta_1;\
\ \rm{ and } $$

$$y_1-z_1+(z_n^2-z_{n-1}^2)(x_{n-1}-y_{n-1})
+z_n(x_{n-1}^2-y_{n-1}^2)=0.$$

\par If $x_{n-1}=y_{n-1}$ , then $x_1=y_1$ and
$\xi_{n-1}=\eta_{n-1}$, so this contribution to
$C_0^t \circ C_0$ is contained in $ \Delta $, the diagonal canonical
relation in $T^*
\R^{n-1} \times T^*\R^{n-1}$.

\par If, on the other hand, $ x_{n-1} \neq y_{n-1}$ ,  then

$$x_i=y_i, \ \ 2 \leq i \leq n-2; $$

$$ \xi_i=\eta_i,  \ \  1 \leq i \leq n-2 ; $$

$$\xi_{n-1}= \frac{(x_{n-1}+y_{n-1})(3x_{n-1}-y_{n-1})}{4}
\xi_1 ;$$

$$\eta_{n-1}= \frac{(x_{n-1}+y_{n-1})(3y_{n-1}-x_{n-1})}{4}
\xi_1; \rm{ and }$$

$$x_1-y_1+ \frac{(x_{n-1}+y_{n-1})^2( x_{n-1}-y_{n-1})}{4}=0,$$
giving a contribution to  $C_0^t \circ C_0$  contained in $ \tilde
{C}_0$, where $ \tilde{C}_0$ is the twisted conormal bundle,

\[\tilde {C}_0= N^* \left\{ x_1-y_1+
\frac{(x_{n-1}+y_{n-1})^2( x_{n-1}-y_{n-1})}{4}=0, x_i-y_i=0,  \ 2
\leq i \leq n-2 \right\}'.\]

\par In conclusion, $ C_0^t \circ C_0 \subseteq \Delta \cup \tilde
{C}_0$, from which
it follows that  $ \tilde{C}_0$ is symmetric, i.e.,
$ \tilde{C}_0^t =  \tilde{C}_0$. It is
easy to see that both projections from $\tilde {C}_0$ have Whitney fold
singularities. Such canonical relations were introduced in \cite{meta}
and called {\it folding canonical relations}; they are
also referred to as {\it two-sided folds} \cite{grse} and we will use
this latter
terminology. One also sees   that
$\tilde {C}_0$ intersects $\Delta$ cleanly in codimension $1$, and
furthermore, $\Delta \cap \tilde {C}_0$ is in fact the fold
surface of $\tilde {C}_0$. As described in $\S 3$, one has a well
defined  class
of distributions associated to the two cleanly intersecting
Lagrangians $(\Delta', \tilde {C}_0')$, namely $I^{p,l} (\mathbb R^{n-1}\times
\mathbb R^{n-1};\Delta',\tilde {C}_0')$, for $p,l \in \R$.  A distribution
is in this class iff it has an oscillatory representation,

\begin{eqnarray*}
u(x,y)&=&\int_{\R^{n-1}} e^{i
\{(x_1-y_1-\frac{(x_{n-1}-y_{n-1})(x_{n-1}+y_{n-1})^2}{4})\xi_1
+(x'''-y''')\cdot\xi'''+(x_{n-1}-y_{n-1})\xi_{n-1}\}} \nonumber\\
& &\hbox{  \ \ \ \ \ \ \ \ \ \ \ \ \ \ \ \ \ \ \ \ \  \ \ \ \ \ \ \ \ 
\ \ \ \ \ }
a(x,y; \xi'';\xi_{n-1}) d\xi
\end{eqnarray*}
where $a$  is a symbol-valued
symbol, satisfying the estimates:

$$|\partial^\alpha_{\xi''}\partial^\beta_{\xi_{n-1}}\partial^\gamma_{x,y}
a(x,y,\xi)| \leq
c(1+|\xi|)^{\tilde{p}-|\alpha|}(1+|\xi_{n-1}|)^{\tilde{l}-\beta},$$
with $\tilde{p}=p+l+\frac12,\tilde{l}=-l-\frac12$.

\par Next we will show that if $F \in I^{m-\frac{1}{4}}(C_0)$ then
$F^*F \in I^{p,l} (\Delta, \tilde {C}_0)$, i.e., $K_{F^*F}\in I^{p,l}(\mathbb
R^{n-1}\times\mathbb R^{n-1};\Delta',\tilde C_0')$, for some
$p,l
\in \R$. We have:
$$Ff(x)=\int e^{
i\{(x''-y'')\cdot\theta''+((x_n^2-x_{n-1}^2)y_{n-1}+ x_n
y_{n-1}^2)\theta_1 \}} a(x,y, \theta) f(y) d\theta dy,$$
where $a\in S^{m+\frac12}$. Thus,
\[ K_{F^*F}(x,y)=\int e^{i\{\phi_0(z,y,\theta'')-\phi_0(z,x,\eta'')\}}
a(z,y,\theta'')\overline{a}(z,x,\eta'') dz d\theta'' d\eta''.
\]
% Hence,
%$F^*F$ has the  phase function
%
%\begin{eqnarray*}
%\phi(x,y,z,\theta, \eta)=
%(z''-y'')\cdot\eta''&+&(z_n^2-z_{n-1}^2)y_{n-1}\eta_1
%+z_ny_{n-1}^2\eta_1\\
%&-&(z''-x'')\theta''
%-(z_n^2-z_{n-1}^2)x_{n-1}\theta_1 -z_nx_{n-1}^2\theta_1.
%\end{eqnarray*}

\par After a stationary phase in the variables $z'',\eta''$, the phase function
becomes:

\begin{eqnarray*}
\phi(x,y, \theta'', z_{n-1}, z_{n})&=&
(x''-y'')\cdot\theta''+(z_n^2-z_{n-1}^2)(y_{n-1}-x_{n-1})\theta_1+
z_n(y_{n-1}^2-x_{n-1}^2)\theta_1\\
&=&(x''-y'')\cdot\theta''+
(y_{n-1}-x_{n-1})[(z_n^2-z_{n-1}^2)+
z_n(x_{n-1}+y_{n-1})]\theta_1
\end{eqnarray*}
and the amplitude becomes $\tilde{a}(x,y,z_{n-1},z_n;\theta'')\in S^{2m+1}$.

\par Following an idea from \cite{gruh4},\cite{fel}, we now  make a
singular change
of variables,
$T: \R^n \rightarrow \R^{n-1},
\ \ T(\theta'', z_{n-1}, z_n)=(\xi'', \xi_{n-1})$, given by:

\begin{eqnarray*}
\xi_i&=&\theta_i,  \ 1 \leq i \leq n-2,\\
\xi_{n-1}&=&-\left((z_n^2-z_{n-1}^2)+
z_n(x_{n-1}+y_{n-1})\right)\theta_1.
\end{eqnarray*}
The kernel of
$F^*F$ can then be rewritten as
$$K_{F^*F}(x,y)=\int_{\R^{n-1}} e^{ i \{ (x''-y'')\cdot\xi''
+(x_{n-1}-y_{n-1})\xi_{n-1}\}}  b(x,y,\xi) d\xi,
$$
where,  using the coarea formula [4, p.249],

\[ b(x,y;\xi) = \int_{ \{((z_n^2-z_{n-1}^2)\theta_1+
z_n(x_{n-1}+y_{n-1})\theta_1)= -\xi_{n-1}\}}  \tilde{a}(x,y, z_{n-1},
z_n;\xi'')
\frac{d\nu}{J_{n-1}}
\]
with $J_{n-1}$ the $(n-1)$-Jacobian of $T$ and $d\nu$ the arc length measure
on\linebreak
$\left\{\left((z_n^2-z_{n-1}^2)+
z_n(x_{n-1}+y_{n-1})\right)\theta_1= -\xi_{n-1}\right\}$.
The  Jacobian is given by
\begin{eqnarray*}
J_{n-1}&=&|\nabla_z
\xi_{n-1}|\\
&=&|( 2z_{n-1} \theta_1, -(2z_n+
x_{n-1}+y_{n-1})\theta_1)|\\
&=& 2
\left(z_{n-1}^2+(z_n+\frac{x_{n-1}+y_{n-1}}{2})^2\right)^{\frac{1}{2}}|\theta_1|,
\end{eqnarray*}
which vanishes to first order at
$z_{n-1}=z_n+\frac{x_{n-1}+y_{n-1}}{2}=0$. Note that, at these points,
$\xi_{n-1}=\frac{(x_{n-1}+y_{n-1})^2}{4}\xi_1$. Thus, $b$ is of order
$2m$ in $\xi''$
and has a conormal singularity  of order $-1$ at
$\{\xi_{n-1}-\frac{(x_{n-1}+y_{n-1})^2}{4}\xi_1=0\}$; thus, it has an
oscillatory
representation
\[ b(x,y,\xi)=\int_{\R}
e^{i\{(\xi_{n-1}-\frac{(x_{n-1}+y_{n-1})^2}{4}\xi_1)\frac{\rho}{\xi_1})\}}
\tilde{b}(x,y; \xi; \rho) d \rho.\]
Finally, the kernel of $F^*F$ becomes

\begin{equation}\label{ker}
K_{F^*F}(x,y)=\int_{\R^n} e^{ i \{ (x''-y'')\cdot\xi''
+(x_{n-1}-y_{n-1})\xi_{n-1} +
(\xi_{n-1}-\frac{(x_{n-1}+y_{n-1})^2}{4}\xi_1)\frac{\rho}{\xi_1}
\}} \tilde{b}(x,y;\xi; \rho) d\xi d\rho
\end{equation}
where one can check that $\tilde{b}$ is a product-type
symbol, $\tilde{b} \in S^{2m, -1}(2n-\nolinebreak2,n-\nolinebreak1,1)$.

We now introduce a conic partition of unity in $(\xi,\rho)$, with supports in
$\{|\rho|\le  2|\xi|\}$ and $\{|\rho|\ge\frac12 |\xi|\}$, resulting in a
decomposition $K_{F^*F}=K^0+K^1$. Letting

$$\psi(x,y; \xi, \rho)= (x''-y'')\cdot\xi'' +(x_{n-1}-y_{n-1})\xi_{n-1}
+(\xi_{n-1}-\frac{(x_{n-1}+y_{n-1})^2}{4}\xi_1)\frac{\rho}{\xi_1},$$
one easily sees that on the region $\{|\rho|\le c |\xi|\}$,  this is
a multi-phase
function for $(\Delta', \tilde{C}_0')$ in the sense of Def.
\ref{old-def3.6}, i.e.,  $\psi_0(x,y; \xi):=\psi_0(x,y, \xi,0)$
parametrizes the
diagonal Lagrangian $\Delta'$ and
$\psi(x,y;(\xi,\rho))$  parametrizes the Lagrangian $\tilde{C}_0'$.
Furthermore, on this region, the amplitude is a symbol-valued symbol,
belonging to
$S^{2m,-1}(\R^{2n-2}\times (\R^{n-1}\setminus 0)\times\R)$.  In view of  Remark
\ref{old-rem3.8},  $K^0\in I^{p,l}(\mathbb R^{n-1}\times\mathbb
R^{n-1};\Delta',
\tilde{C}_0')$,  for some $p,l\in\R$.
The orders may also be found
by applying Remark \ref{old-rem3.8} to (\ref{ker}). We see that
\[
p=(2m-1) +\frac{(n-1)+1}2-\frac{2n-2}4=2m-\frac12
\]
and
\[
l=-(-1)-\frac12=\frac12.
\]
Hence, this contribution to $F^*F$ is in
$I^{2m-\frac{1}{2}, \frac{1}{2}}(\Delta,\tilde{C}_0)
:=I^{2m-\frac{1}{2}, \frac{1}{2}}(\R^{n-1},\R^{n-1};\Delta,\tilde{C}_0)$,
with $\tilde{C}_0$  a two-sided fold.

Next, we show  that  $K^1\in I^{2m-\frac{1}{2}}(\tilde {C}_0)$.
First, let $s=\xi_{n-1}/\xi_1$. Then we can express
$K^1(x,y)= \int L(x,y,s) ds$ ,
where $L$ has the phase function
$$\Psi(x,y,s;\rho;
\xi'')=(x''-y'')\cdot\xi'' +(x_{n-1}-y_{n-1})\xi_1 s
+(s-\frac{(x_{n-1}+y_{n-1})^2}{4})\rho$$
and amplitude $c\in S^{-1,2m+1}(\R^{2n-1}\times (\R\setminus
0)\times\R^{n-1})$.
Note that
$\Psi_0(x,y,s;\rho):=\Psi(x,y,s;\rho;0)$ parametrizes
\[
\Gamma_0:= N^*\{
s-\frac{(x_{n-1}+y_{n-1})^2}{4}=0 \}\subset T^*\R^{2n-1}\setminus 0
\]
and
$\Psi_1(x,y,s;(\rho,\xi'')):=\Psi(x,y,s;(\rho,\xi''))$ parametrizes
\[
\Gamma_1:= N^* \left\{ x_i-y_i=0,  2 \leq i\leq n-2,
\ x_1-y_1+s(x_{n-1}-y_{n-1})=
s-\frac{(x_{n-1}+y_{n-1})^2}{4}=0 \right\}.
\]
Thus, $\Psi$ is
a multi-phase function in the sense of Def $3.5$.
Introduce a cutoff
function, $\chi \in C_0^{\infty}(R^n), \chi=1 $
near $0$,   set
$$L^0(x,y,s)=\int e^{i \Psi(x,y,s;\rho,
\xi'')}\chi(\frac{|\xi''|}{c|\rho|})
c(x,y;\rho;\xi'') d\xi'' d\rho,$$
and let $L^1= L- L^0$. We have $L^0 \in
I^{2m-\frac{1}{4},\frac{1}{2}} (\R^{2n-1};\Gamma_0, \Gamma_1)$, $L^1 \in
I^{2m-\frac{1}{4}} (\R^{2n-1};\Gamma_1)$ and $WF(L^1)$ is contained
in the complement
of a neighborhood of $\Gamma_0$.

\par Since  $K^2$ is
simply $\pi_{*}L$, the pushforward of $L$ by $\pi$, the projection
$\pi(x,y,s)=(x,y)$, which  is a submersion, we compute the pushforwards of the
Lagrangians $\Gamma_0,
\Gamma_1$.  It follows from standard results about pushforwards \cite{ho}
that, for
$u\in\mathcal E'\left(\R^{2n-1}\right)$,
\begin{eqnarray*}
WF(\pi_{*}u) \subseteq \{ (x,y; \xi, \eta)&|&
\exists (\hat{x},
\hat{y}, s, \hat\xi, \hat\eta, \sigma) \in WF(u) \ \ s.t.\\
&{ }&(x,y) =
\pi(\hat{x},
\hat{y},s), (\hat\xi, \hat\eta, \sigma)= d \pi^t (\xi,
\eta) \}.
\end{eqnarray*}
Using $d \pi^t (\xi,\eta)=(\xi,\eta,0)$, it is then  an easy
calculation that the
push forward of
$\Gamma_0$, and indeed any neighborhood of $\Gamma_0$ on which
$\sigma\ne0$, is empty,
so that
$\pi_*(L^0)\in C^\infty$, while
   $\pi_*(\Gamma_1)\subseteq\tilde{C}_0'$.
In fact, $\pi_{*}$ is a FIO of order $- \frac{1}{4}$, and the application of
$\pi_*$ to $I^\cdot(\Gamma_1)$ is covered by the transverse intersection
calculus. Hence
$\pi_{*}: I^{r}(\Gamma_1)
\rightarrow
I^{r-\frac{1}{4}}(\tilde{C}_0')$, so that
\[
K^1=\pi_{*}(L^1) \mod C^\infty\in I^{2m-\frac{1}{2}} (\tilde{C}_0),
\]
and thus
$K_{F^*F}=K^0+K^1 \in I^{2m-\frac{1}{2},\frac12}(\Delta',\tilde{C}_0')$.
\vskip.2in
In conclusion, we have shown that $F^*F\in 
I^{2m-\frac{1}{2},\frac12}(\Delta,\tilde{C}_0)$.
For the single source data acquisition
geometry,  Nolan\cite{no} and Felea\cite{fel}  showed that  $F^*F \in
I^{2m,0} (\Delta, \tilde {C})$, with $\tilde{C}$ a two-sided fold. In
that situation the artifact, i.e., the part of $F^*F$ on
$\tilde{C} \setminus \Delta$,
has, by Remark \ref{old-rem3.4}, the same strength as on
$\Delta\setminus\tilde{C}$.
However, it follows from what we have shown here that for the
microlocal model $C_0$
of  the marine geometry, the artifact  is $\frac{1}{2}$ order smoother then the
pseudodifferential operator part. In the next section, we show
this in full generality for FIOs associated with folded cross caps.

\section{Composition calculus in the general  case}

\par By a  well known result of Melrose and
Taylor\cite{meta},
any two-sided fold can be conjugated, via canonical transformations
on the left and
right, to a normal form. Unfortunately,  even if  folded cross caps  could be
conjugated to a  normal form, such a result would presumably be difficult to
prove. However, we will be able to avoid this problem by finding a
weak normal form
for the canonical relation and thus for a phase function
parametrizing it, adapting a
method originally developed by Greenleaf and Uhlmann \cite{gruh3,gruh4} for
some canonical
relations arising in integral geometry, for which
$\pi_R$ belongs to a class containing the submersions with folds and
$\pi_L$ is maximally degenerate. This will be sufficient
for establishing the  composition calculus for general folded cross
cap canonical
relations.
\vskip.2in

\begin{theorem}\label{main} {\em If} $F \in I^{m-\frac{1}{4}}(X,Y;C)$ {\em is
properly  supported
and} $C$  {\em is a folded cross cap canonical relation, then}
$F^*F \in I^{2m- \frac{1}{2},\frac{1}{2}}(Y,Y;\Delta, \tilde{C})$.
   {\em Furthermore, }  $\tilde
{ C}$ {\em is a symmetric, two-sided fold in} $\left(T^*Y \setminus
0\right) \times
\left(T^*Y \setminus 0\right)$, {\rm and }
$\Delta\cap\tilde C$ {\rm equals the fold surface in } $\tilde C$.
\end{theorem}
\vskip.2in
\begin{remark}
It was shown in \cite{nosy} that, for the marine seismic data 
geometry in three spatial dimensions, the linearized scattering map 
$F$ belongs to $I^{1-\frac14}(\Sigma_{r,s}\times (0,T),\R^3; C)$. 
Thus, in the presence of fold caustics and under the additional 
nonradiality assumptions described below
Def. \ref{def-fcc},  it follows from Theorem \ref{main} that the 
normal operator $F^*F$ lies in 
$I^{\frac32.\frac12}(\R^3,\R^3;\Delta,\tilde{C})$, with $\tilde{C}$ 
as above.
\end{remark}
\vskip.2in

Before establishing a weak normal form for a general folded cross
cap, we first need
to find separate weak normal forms for each of the two projections,
$\pi_R$ and $\pi_L$.

\begin{proposition}\label{sympl-swf}
Let $f:V \longrightarrow W$ be a smooth, conic map, with $V$ a
smooth, conic manifold
of dimension $2n-1$ and $(W,\omega_W)$ a conic symplectic manifold of dimension
$2n-2$. Assume that $f$ is a submersion with folds at $v_0\in V$ and $f(V)$ is
nonradial in $W$. Then, there exist local homogeneous coordinates
$(s'',s_{n-1},\sigma'',\sigma''')\in\R^{n-2}\times\R\times(\R^{n-2}\setminus
0)\times\R^2$ on $V$ near $v_0$, and local canonical coordinates
$(y',\eta')\in T^*\R^{n-1}\setminus 0$ on $W$ near $w_0:=f(v_0)$, such that
$v_0=(0,0,e_1^*,0),\quad w_0=(0;e_1^*)$, and
\begin{equation}\label{swfnf}
f(v)=f(s'',s_{n-1},\sigma'',\sigma''')
=\left(s'',s_{n-1};\sigma'',A^{s,\sigma}(\sigma''',\sigma''')\right),
\end{equation}
where $A^\cdot\in C^\infty(\R^{2n-1},{S^2\R^2}^*)$ is homogeneous of
degree -1 in
$\cdot$ and  takes values in the nonsingular quadratic forms of the
same signature as
Hess $f(v_0)$.
\end{proposition}

\noindent\emph{Proof.}
In the terminology of \cite{gruh4}, a submersion with folds has {\it
clean folds of multiplicity one}.  Prop.\ref{sympl-swf} is then a
special case of
\cite[Lem. 3.A.1]{gruh4}, with slight change of notation, and  $(N,m,n)$ in
\cite{gruh4} being $(n-2,1,2)$ here. $\square$

\vskip.2in
Applying (\ref{swfnf}) to $\pi_R:C\lra T^*Y\setminus 0$ and writing
$\sigma'''=(\sigma_{n-1},\sigma_n)$, we see that
${\cal S}_1^C={\cal S}_1(\pi_R)=\{\sigma_{n-1}=\sigma_n=0\}$ and $\pi_R({\cal
S}_1^C)=\{\eta_{n-1}=0\}$. Furthermore,
\begin{eqnarray}\label{omegc}
\omega_C &=& \pi_R^*(\sum_{j=1}^{n-1} d\eta_j\wedge dy_j)\nonumber\\
&=& \sum_{j=1}^{n-2} d\sigma_j\wedge ds_j + d\left(A^{s,\sigma)}
(\sigma''',\sigma''')\right)\wedge ds_{n-1}.
\end{eqnarray}
Hence,
\begin{equation}\label{ker-swf}
Ker \left(\omega_C|_{T{\cal S}_1}\right)=\mathbb
R\cdot\left(\frac{\partial}{\partial
s_{n-1}}+\dots\right).
\end{equation}
Also, as in the case of the model canonical relation, $C_0$, we see that
$\omega_C$ has rank $2n-2$ on $C\setminus {\cal S}_1^C$, and has rank
$2n-4$ both at
${\cal S}_1^C$ and on ${\cal S}_1^C$, i.e., restricted to $T{\cal S}_1^C$.
Thus, since $\omega_C=\pi_L^*\omega_{T^*X}$ as well, the image $\pi_L({\cal
S}_1)\subset T^*X\setminus 0$, which is smooth, conic, nonradial and
codimension 3,
must also be  maximally noninvolutive. Finally, we note that, as a general fact
about canonical relations, for all $c_0=(x_0,\xi_0,y_0,\eta_0)\in C$,
the subspace
$d\pi_L(T_{c_0}C)\le T_{(x_0,\xi_0)}\left(T^*X\right)$ is involutive;
in particular,
for $c_0\in{\cal S}_1$, $d\pi_L(T_{c_0}C)$ is a codimension 2,
involutive subspace.
We
are thus led to establishing a (very) weak normal form for cross cap maps into
symplectic manifolds reflecting these extra conditions. Note that
this would apply to
a more general class of maps than cross caps, since at this point we
are only using
information concerning the first derivatives.
\vskip.2in
\begin{proposition}\label{sympl-cc}
Let $g:V\lra U$ be a smooth conic map, with $V$ a
smooth, conic manifold of dimension $2n-1$ and $(U,\omega_U)$ a conic
symplectic
manifold of dimension $2n$. Assume that $g$ has a cross cap singularity, with
critical set ${\cal S}_1$. Let $v_0\in{\cal S}_1, u_0=g(v_0)$. Assume  that
$dg(T_{v_0}V)\le T_{u_0}(U)$ is involutive  and $dg\left(T_{v_0}{\cal
S}_1\right)=T_{u_0}\left(g({\cal S}_1)\right)\le T_{u_0}U$ is maximally
noninvolutive for all $v_0\in{\cal S}_1$. Then, there exist local homogeneous
coordinates
\[
(t,\tau)=(t'',t_n,\tau'',\tau''')\in\R^{n-2}\times\R\times(\R^{n-2}
\setminus 0)\times\R^2
\]
on $V$ near $v_0$, and local canonical coordinates $(x,\xi)\in
T^*\R^n\setminus 0$ on
$U$ near $u_0:=g(v_0)$, such that $v_0=(0,0,e_1^*,0)$,
$u_0=(0;e_1^*)$, and, writing
$\tau'''=(\tau_{n-1},\tau_n)$,
\begin{equation}\label{ccnf}
g(t,\tau)=\left(t'',g_{n-1}(t,\tau),t_n;\tau'',
\gamma_{n-1}(t,\tau),\gamma_n(t,\tau)\right);
\end{equation}

\begin{equation}\label{gnonz}
\frac{\partial g_{n-1}}{\partial \tau_{n-1}}\ne 0;\quad {\text and}
\end{equation}

\begin{equation}\label{derivz}
\gamma_i|_{\tau'''=0}=\frac{\partial
\gamma_i}{\partial\tau_j}|_{\tau'''=0}=0,\quad
n-1\le i,j\le n.\end{equation}
\end{proposition}
\vskip.2in

Note that, with respect to these coordinates,
$\so\subseteq\{\tau_{n-1}=\tau_n=0\}$,
$\Sigma^{2n-3}:=g(\so)\subseteq\{x_{n-1}=\xi_{n-1}=\xi_n=0\}$ and

\begin{equation}\label{ker-cc}
Ker\left(\omega_C|_{T\so}\right)=\R\cdot\left(\frac{\partial}{\partial
\tau_{n-1}}+\dots\right).
\end{equation}
\vskip.2in

\noindent\emph{Proof.}  By assumption, $\Sigma^{2n-3}\subset U$ is nonradial,
codimension 3 and maximally nonivolutive. By an application of
Darboux's Theorem,
there exist (\cite[Thm. 21.2.4]{ho-book}) local canonical coordinates
$(x,\xi)\in
T^*\R^n\setminus 0$ such that
$\Sigma^{2n-3}\subseteq\{x_{n-1}=\xi_{n-1}=\xi_n=0\}$
and
$u_0=(0,e_1^*)$. Letting $t_j=g^*(x_j),\quad \tau_j=g^*(\xi_j)$ for
$1\le j\le n-2$,
and $t_n=g^*(x_n)$, these functions have linearly independent
gradients at $v_0$. If
$\tau_{n-1},\tau_n$ are any two defining functions for $\so$,
homogeneous  of degree
1, then $(t,\tau):=(t'',t_n,\tau'',\tau_{n-1},\tau_n)$ are local homogeneous
coordinates on $V$, with $\so\subseteq\{\tau_{n-1}=\tau_n=0\}$. With
respect to these
coordinates,
\[
g(t,\tau)=\left(t'',x_{n-1}(t,\tau),t_n;\tau'',\xi_{n-1}(t,\tau),\xi_n(t,\tau)\right)
\]
and $\rank dg(v) = (2n-3)+\rank \mathbf D(v)$, where
\[
\mathbf D(v)= \left(\begin{array}{cc}
\frac{\partial x_{n-1}}{\partial \tau_{n-1}} &
\frac{\partial x_{n-1}}{\partial \tau_n} \\
{ } & { } \\
\frac{\partial \xi_{n-1}}{\partial \tau_{n-1}} &
\frac{\partial \xi_{n-1}}{\partial \tau_n} \\
{ } & { } \\
\frac{\partial \xi_{n}}{\partial \tau_{n-1}} &
\frac{\partial \xi_{n}}{\partial \tau_n}
\end{array}\right) = \left( D_{n-1},D_n\right),\quad D_{n-1},D_n\in T_uU.
\]
Since $g$ is a cross cap, $\rank \mathbf D(v)=1,\forall v\in\so$, so that by
interchanging $\tau_{n-1}$ and $\tau_n$, if necessary, we can assume
that $D_{n-1}\ne
0$ and $D_n\in\R\cdot D_{n-1}$ for $v\in\so$ near $v_0$. Furthermore,
rotating in the
$x_{n-1},\xi_{n-1}$ plane, if necessary, we can likewise assume that
$\frac{\partial
x_{n-1}}{\partial \tau_{n-1}}(v_0)\ne 0$ and
$\frac{\partial \xi_{n-1}}{\partial \tau_{n-1}}(v_0)= 0$, so that
\[
\left|\frac{\partial x_{n-1}}{\partial \tau_{n-1}}(v)\right| >>
\left|\frac{\partial \xi_{n-1}}{\partial \tau_{n-1}}(v)\right|,\quad\forall v
{\textrm{ near }} v_0.
\]
Let
\[
\Pi(v)=dg(T_vV)=span\left\{
\{\frac{\partial}{\partial x_j}\}_{j=1}^{n-2},\frac{\partial}{\partial x_n},
\{\frac{\partial}{\partial \xi_j}\}_{j=1}^{n-2}, D_{n-1}
\right\};
\]
by assumption, $\mathbf\Pi(v)\le T_{g(v)}U$ is a codimension 2,
involutive subspace.
We have
\[
D_{n-1}(v)=a(v)\frac{\partial}{\partial x_{n-1}}
+ b(v) \frac{\partial}{\partial \xi_{n-1}}
+ c(v) \frac{\partial}{\partial \xi_{n}},
\]
with $|a|>>|b|$ near $v_0$. A simple calculation shows that a vector
\linebreak$\mathbf X=\sum_{j=1}^n \alpha_j \frac{\partial}{\partial x_{j}}
+ \beta_j \frac{\partial}{\partial \xi_{j}}\in T_{g(v)}U$ is in the symplectic
annihilator
$dg(T_vV)^\omega$ iff
\[
\alpha_j=\beta_j=0,\quad 1\le j\le n-2, \beta_n=0,\quad
a\beta_{n-1}-b\alpha_{n-1}-c\alpha_n=0.
\]
Thus, $\left\{b\frac{\partial}{\partial x_{n-1}} + a\frac{\partial}{\partial
\xi_{n-1}}, c\frac{\partial}{\partial x_{n-1}}+ a
\frac{\partial}{\partial \xi_n}
\right\}$ forms a basis for $dg(T_vV)^\omega$. Since $dg(T_vV)$ is involutive
iff $dg(T_vV)^\omega\le dg(T_vV)$, we must have $b(v)=c(v)=0$. Thus,
$\frac{\partial \xi_i}{\partial \tau_j}|_{\tau'''=0}=0,\quad n-1\le i,i\le n$.
Relabeling, this finishes the proof of Prop.\ref{sympl-cc}. $\square$

\vskip.2in

Now, let $C\subset (T^*X\setminus 0)\times (T^*Y\setminus 0)$ be a
folded cross cap
canonical relation, and apply both Prop.\ref{sympl-swf}  to
$\pi_R:C\lra T^*Y\setminus 0$ and Prop.\ref{sympl-cc} to $\pi_L:C\lra
T^*X\setminus
0$ near $c_0\in C$. As in \cite{gruh3,gruh4}, we now attempt to
reconcile the two
resulting coordinate systems on $C$, $(s'',s_{n-1},\sigma'',\sigma''')$ and
$(t'',t_n,\tau'',\tau''')$. On $T\so$ (which is the same for both projections),
$\omega_C$ has rank $2n-4$. By (\ref{ker-swf}) and (\ref{ker-cc}), since the
hypersurface $\{\tau_{n-1}=0\}$ is transverse to $\Ker(\omega_C)$, it
must be locally
expressible as $\{s_{n-1}=\tilde{F}(s'',\sigma'')\}$ for some smooth
$\tilde{F}$. Now
set
$F(s,\sigma)=-A^{s,\sigma}(\sigma''',\sigma''')\tilde{F}(s'',\sigma'')$;
then
$F=\pi_R^*f$, where $f(y,\eta)=-\eta_{n-1}\tilde{F}(y'',\eta'')$. Using
(\ref{omegc}), we can find a vector field $\mathbf Y$ on $C$ which satisfies
\[
i(\mathbf Y)\omega_C= dF =
-\tilde{F}(s'',\sigma'')\frac{\partial}{\partial s_{n-1}}
+ O(|\sigma'''|^2)\cdot (ds'',d\sigma'').
\]
Thus, $exp(H_f)$ is a canonical transformation of $T^*Y\setminus 0$, while
$exp(\mathbf Y)$ is an $\omega_C$-morphism of $C$, mapping $\{\tau_{n-1}=0,
\tau'''=0\}$ into $\{s_{n-1}=0,\sigma'''=0\}$. Applying these simultaneously on
$T^*Y$ and $C$, respectively, we see that one can assume that
$L:=\{\tau_{n-1}=0,\tau'''=0\}=\{s_{n-1}=0,\sigma'''=0\}$. Restricted to this
$(2n-4)$-dimensional submanifold, $\omega_C$ is symplectic, so by Darboux there
exists a canonical transformation
$\Phi(y'',\eta'')=(\Phi_{y''},\Phi_{\eta''})$ of
$T^*\R^{n-2}\setminus 0$ such that
\[
\Phi^*( s''|_L)= t''|_L\quad\textrm{ and } \Phi^*(\sigma''|_L)=\tau''|_L.
\]
Extend $\Phi$ to  $\tilde{\Phi}:T^*Y\setminus 0\lra T^*Y\setminus 0$ by
\[
\tilde{\Phi}(y'',y_{n-1};\eta'',\eta_{n-1})=(\Phi_{y''}(y'',\eta''),
y_{n-1};\Phi_{\eta''}(y'',\eta''),\eta_{n-1})
\]
and compose $C$ on the right with the graph of $\tilde{\Phi}$. Then,
\[
\pi_L^*(x'')=t''=s'', \pi_L^*(x_{n-1}),
\pi_L^*(x_n)=t_{n},\pi_R^*(y_{n-1})=s_{n-1}\textrm{ and
}\pi_R^*(\eta'')=\sigma''
\]
form coordinates on $C$ near $c_0$.

\par Thus,  we have so far shown that if $C$ is a folded cross cap,
we may assume that
$(x, y_{n-1},
\eta'')$ form (micro)local coordinates on
$C$. Therefore  there is  a generating function
$S(x, y_{n-1},\theta'')$ for
$C$, where
$S$ is $C^{\infty}$ and homogeneous of degree $1$ in $\theta''$
(\cite[Thm.21.2.18]{ho-book}). Hence, $C$ can be parametrized as
\begin{equation}\label{i}
C=\left\{ (x, d_x S; d_{\theta''} S,
y_{n-1}, -\theta'', -d_{y_{n-1}} S) \right\}
\end{equation}
and the phase function
$\chi(x,y,\theta'')=S(x,y_{n-1},\theta'')-y'' \cdot \theta''$ is a
parametrization for the Lagrangian $C'$. We now show that the
properties of $\pi_L$
and $\pi_R$ impose several conditions  on $S$ and its derivatives,
forcing the phase
function to be very similar to the model phase $\phi_0$ considered in \S5.

\par To prepare the canonical relation $C$, we first replicate (\ref{ii}) and
(\ref{iii}).
Since $\pi_R({\cal S}_1^C)$ is a nonradial hypersurface in $T^*Y \setminus
0$, by Darboux's  theorem,  microlocally there
is a canonical transformation from $T^*Y \setminus 0$ to $T^*\R^{n-1}
\setminus 0$ taking  $\pi_R({\cal S}_1^C)$   to
$\pi_R({\cal S}_1^{C_0})$, given by (\ref{ii}). Similarly,
$\pi_L({\cal S}_1^C)$ is a codimension three submanifold of $T^*X\setminus 0$,
which is maximally noninvolutive in the sense that $\omega_{T^*X}
|_{\pi_L({\cal
S}_1^C)}$ has rank $2n-4$. Thus, there exist  defining functions $p_1,
p_2, p_3$ for $\pi_L({\cal S}_1^C)$ with $\{p_1, p_2 \}=1, \{p_1,p_3 \}=0,
\{p_2,p_3 \}=0$. By Darboux's theorem, we can find a
canonical transformation from $T^*X \setminus 0$ to $T^*\R^n
\setminus 0$ mapping $\pi_L({\cal S}_1^C)$ to $\pi_L({\cal S}_1^{C_0})$
   given by (\ref{iii}).

\par Using (\ref{i}), $\pi_R(x, y_{n-1}, \theta'')=(d_\theta'' S, y_{n-1},
-\theta'', -d_{y_{n-1}}
S)$, so (\ref{ii}) implies

\begin{equation}\label{old-one}
-d_{y_{n-1}} S|_{\{x_{n-1}=0=x_n+y_{n-1}\}}=x_n^2\theta_1
\end{equation}

\par From (\ref{i}), $\pi_L (x, y_{n-1}, \theta'')=(x, d_x S)$ and
(\ref{iii}) implies

\begin{equation}\label{old-two}
d_{x_{n-1}} S|_{\{x_{n-1}=0=x_n+y_{n-1}\}} =0
\end{equation}
and

\begin{equation}\label{old-three}
d_{x_n} S |_{\{x_{n-1}=0=x_n+y_{n-1}\}}=-x_n^2 d_{x_1} S |_{\{
x_{n-1}=0=x_n+y_{n-1} \}}
\end{equation}

\par Relation (\ref{old-two}) means that

\begin{equation}\label{old-four}
S(x, y_{n-1}, \theta'')= S_0(x'',\theta'')+x_{n-1}^2
S_1(x, y_{n-1}, \theta'')+ (x_n+ y_{n-1}) S_2(x, y_{n-1}, \theta'').
\end{equation}

\par  From  (\ref{old-one}) we obtain

\begin{equation}\label{old-five}
S_2(x, y_{n-1}, \theta'')|_{ \{x_{n-1}=0=x_n+y_{n-1}
\} }=-x_n^2\theta_1 ,
\end{equation}
and hence

\begin{equation}\label{old-three}
S_2(x, y_{n-1}, \theta'')= - x_n^2\theta_1 + x_{n-1}
S_3(x, y_{n-1}, \theta'') + (x_n + y_{n-1}) S_4(x, y_{n-1}, \theta'').
\end{equation}

\par Identity (\ref{old-three}) implies that
\begin{equation}\label{old-six}
S_2(x, y_{n-1},
\theta'')|_{\{x_{n-1}=0=x_n+y_{n-1}\}}=- x_n^2 d_{x_1}
S|_{\{x_{n-1}=0=x_n+y_{n-1} \}}.
\end{equation}

\par Now, (\ref{old-five}) and (\ref{old-six}) imply  $ d_{x_1}
S|_{\{x_{n-1}=0=x_n+y_{n-1}\}}= \theta_1$, and from (\ref{old-four}) we have
that $d_{x_1} S_0 =\theta_1$, so $S_0=x_1 \theta_1  +
\tilde{S}_0(x''', \theta''')$.  At this point, our analysis shows
that the generating
function has the form

\begin{eqnarray}\nonumber
S(x, y_{n-1}, \theta'')= x_1\theta_1&+&\tilde{S}_0 (x''', \theta''')+
x_{n-1}^2 S_1(x, y_{n-1}, \theta'') \\
&+& (x_n+y_{n-1})[-x_n^2\theta_1
+x_{n-1} S_3(x, y_{n-1}, \theta'')\label{old-seven}\\
& { } & { }  +( x_n+ y_{n-1}) S_4(x, y_{n-1},
\theta'') ] .\nonumber
\end{eqnarray}

\par Next we consider the differentials $d\pi_R$ and $d\pi_L$. One computes

\[
{d \pi_R}= \left(\begin{array}{ccccc}

d_{x''\theta''}S & d_{x_{n-1}\theta''}S & d_{x_n \theta''}S &
d_{y_{n-1}\theta''}S & d_{\theta''\theta''}S\\
0 & 0 & 0 & 1 & 0\\
0 & 0 & 0 & 0 & -I_{n-2}\\
-d_{x''y_{n-1}}S & -d_{x_{n-1}y_{n-1}}S &  -d_{x_n y_{n-1}}S &
-d_{y_{n-1}y_{n-1}}S & -d_{\theta''y_{n-1}}S
\end{array} \right)
\]

\par Evaluating at ${\cal S}_1^C$, by (\ref{old-seven}), this becomes \\

\[{d \pi_R|_{\so^C}}= \left(\begin{array}{ccccc}
d_{x''\theta''}S_0 & 0 & 0 & 0 & d_{\theta''\theta''}S_0\\
0 & 0 & 0 & 1 & 0\\
0 & 0 & 0 & 0 & -I_{n-2}\\
0 & S_3 &  -2x_n\theta_1 + 2S_4 & 2S_4 & 0
\end{array} \right)\]

\par By assumption, $d\pi_R$ has rank $2n-3$ at ${\cal S}_1^C$.
Thus,  $S_0(x'',\theta'')$ is
nondegenerate,   $S_3|_{ \{x_{n-1}=0=x_n+y_{n-1} \}}$=0 and
$S_4|_{ \{ x_{n-1}=0=x_n+y_{n-1} \}}=x_n \theta_1$. We thus have

\begin{equation}\label{old-eight}
S_3(x, y_{n-1}, \theta'') =x_{n-1}S_5(x, y_{n-1},
\theta'')+(x_n+y_{n-1})S_6(x, y_{n-1}, \theta''),
\end{equation}
and

\begin{equation}\label{old-nine}
S_4(x, y_{n-1}, \theta'')= x_n \theta_1 +
x_{n-1}S_7(x, y_{n-1}, \theta'')+(x_n+y_{n-1})S_8(x, y_{n-1}, \theta'').
\end{equation}

\par Similarly,

\[{d \pi_L|_{\so^C}}= \left(\begin{array}{ccccc}
I_{n-2} & 0 & 0 & 0 & 0 \\
0 & 1 & 0 & 0 & 0  \\
0 & 0 & 1 & 0 & 0  \\
d_{x''x''}S & d_{x_{n_1}x''}S & d_{x_nx''}S & d_{y_{n-1}x''}S &
d_{\theta''x''}S   \\
d_{x''x_{n-1}}S & d_{x_{n-1}x_{n-1}}S & d_{x_nx_{n-1}}S &
d_{y_{n-1}x_{n-1}}S &  d_{\theta''x_{n-1}}S\\
d_{x''x_n}S & d_{x_{n-1}x_n}S & d_{x_nx_n}S & d_{y_{n-1}x_n}S &
d_{\theta''x_n}S
\end{array} \right).\]

\par Using relation (\ref{old-seven}) again, we have

\[{d \pi_L|_{\so^C}}= \left(\begin{array}{ccccc}
I_{n-2} & 0 & 0 & 0 & 0 \\
0 & 1 & 0 & 0 & 0  \\
0 & 0 & 1 & 0 & 0  \\
d_{x''x''}S_0& 0 & 0 & 0 & d_{\theta''x''}S_0   \\
0 & 2S_1 & N & S_3 &  0\\
0 & N & 2S_4-2x_n \theta_1 & 2S_4-2x_n\theta_1 & 0
\end{array} \right),\]
where $ N= S_3+\partial_{x_{n-1}}S_4$.

\par By the cross cap condition,  the rank of
$d \pi_L|_{\{x_{n-1}=0=x_n+y_{n-1} \}}$ is
$2n-2$. Thus, $S_3|_{\{x_{n-1}=0=x_n+y_{n-1} \}}=0$ and
$S_4|_{\{x_{n-1}=0=x_n+y_{n-1} \}}= x_n \theta_1$, and  we
see that  relations (\ref{old-eight}) and (\ref{old-nine}) follow 
from the analysis of $\pi_L$ as well as
$\pi_R$.
Putting together all
the previous relations, the generating function  becomes:

\begin{eqnarray}\label{old-ten}
S(x, y_{n-1}, \theta'')=
x_1\theta_1&+&\tilde{S}_0 (x''', \theta''')+ x_{n-1}^2 S_1(x,
y_{n-1}, \theta'')
\nonumber\\
& +&  (x_n+y_{n-1})x_ny_{n-1}\theta_1\nonumber\\
&+&(x_n+y_{n-1})[ x_{n-1}^2 S_5(x, y_{n-1}, \theta'')\\
&{ }& \ \  \ \ \ \ \ \ \ \ +x_{n-1}(x_n+ y_{n-1}) S_6(x, y_{n-1}, 
\theta'')\nonumber\\
&{ }& \ \ \ \ \ \ \ \ \ \ +(x_n+y_{n-1})^2 S_7(x,y_{n-1}, \theta'')].\nonumber
\end{eqnarray}

\par Since $S_0$ is nondegenerate, $C \cap \{x_{n-1}=x_n=y_{n-1}=0
\}$ is  the graph of a canonical transformation on
$T^*R^{n-2}$; we may thus  assume that  $C \cap
\{x_{n-1}=x_n=y_{n-1}=0 \}=\left\{ (x'', 0, \theta'', 0; x'', 0,
\theta'',0) \right\}$, so that $S_0(x'', \theta'') =x'' \cdot \theta''$,
and (\ref{old-ten}) becomes:

\begin{eqnarray*}
S(x, y_{n-1}, \theta'')=
x''\cdot\theta''&+& x_{n-1}^2 S_1(x, y_{n-1}, \theta'')+
(x_n+y_{n-1})x_ny_{n-1}\theta_1 \\
&+&(x_n+y_{n-1})[ x_{n-1}^2
S_5(x, y_{n-1}, \theta'')\\
&{ }& \ \ \ \ \ \ \ \ \ \ x_{n-1}(x_n+ y_{n-1}) S_6(x,
y_{n-1}, \theta'')\\
&{ }& \ \ \ \ \ \ \ \ \ \ (x_n+y_{n-1})^2 S_7(x,y_{n-1}, \theta'')]
\end{eqnarray*}

\par  Letting

\begin{eqnarray*}
M(x,y_{n-1},\theta'')=x_{n-1}^2S_5(x,y_{n-1},\theta'')
+x_{n-1}(x_n+y_{n-1})S_6(x,y_{n-1},\theta'')\\
+(x_n
+ y_{n-1})^2S_7(x,y_{n-1},\theta''),
\end{eqnarray*}
we now see that the Lagrangian $C'$ is  parametrized by the phase function

\begin{eqnarray}\label{old-eleven}
\chi(x,y,\theta'')&=&S(x,y_{n-1},\theta'')-y''\cdot\theta''\nonumber\\
&=& (x''-y'')\cdot\theta''+ \left(x_n^2 y_{n-1}
+x_ny_{n-1}^2\right)\theta_1\\
& { } & \ \ \ \  + x_{n-1}^2 S_1(x,y_{n-1}, \theta'' )
+(x_n+y_{n-1})M(x,y_{n-1},\theta'').\nonumber
\end{eqnarray}

\par The right projection becomes:
$$
\pi_R\big(x'', y_{n-1}, \theta'', x_{n-1}, x_n)=(x''+ x_{n-1}^2
\partial_{\theta''}S_1 + (x_n+y_{n-1})\partial_{\theta''}M,
y_{n-1}, \theta'',$$ $$ (x_n^2+2x_n y_{n-1})\theta_1+x_{n-1}^2
\partial_{y_{n-1}}S_1 + M + (x_n+y_{n-1})\partial_{y_{n-1}}M \big)
$$

\par Since $\pi_R$ is a submersion with folds,
$Hess(d\pi_R):Ker(d\pi_R)\lra Coker(d\pi_R)$,
its Hessian
at $\so^C$, is
nonsingular. At $\so^C$, we have $Ker(d\pi_R)=span\{\frac{\p}{\p
x_{n-1}},\frac{\p}{\p x_n}\}$ and

$$ Hess(\pi_R)= \left(\begin{array}{cc}
2 \partial_{y_{n-1}}S_1|_{S_1^C} + 2 S_5|_{S_1^C} & 2S_6|_{S_1^C} \\
2S_6|_{S_1^C} & 2 \theta_1 + 6 S_7|_{S_1^C}
\end{array} \right)$$

\par At a fixed point $c_0\in\so^C$, by a rotation in $x_{n-1},x_n$,
we can  diagonalize this to obtain $S_6=0$ at $c_0$. Similarly,
$2\theta_1+6S_7$ represents the
$d\xi_n$ component of $Hess(\pi_L)$; by a canonical transformation,
we can assume that this equals $2\theta_1$ at $c_0$, so that $S_7$
vanishes at $c_0$. Thus,
$S_6$ and $S_7$ are small near $c_0$; note
also that, by the nondegeneracy of $Hess(\pi_R)$,
$\partial_{y_{n-1}}S_1 + 2 S_5 \neq 0$ near $c_0$.
\medskip

  We are now ready to prove Theorem \ref{main}. Composing $F$ on the
left and right with elliptic FIOs of order $0$ associated with all
of the canonical transformations of $T^*X \setminus 0$ and $T^*Y
\setminus 0$ used above, we can assume that the Schwartz kernel of
$F$ is represented by an oscillatory integral with the phase
function given by (\ref{old-eleven}) and an amplitude $a\in S^{m+\frac12}$.

\par Let $\tilde {\chi}$ be the phase function of $F^*F$:

\begin{eqnarray*}
\tilde {\chi}&=&\chi(z,y,\eta'')-\chi(z,x,\xi'')\\
&=&(z''-y'')\cdot\eta'' +z_n^2y_{n-1}\eta_1+ z_ny_{n-1}^2\eta_1
+ z_{n-1}^2 S_1(z, y_{n-1}, \eta'') \\
& { } &  \ \ \ \ +(z_n+y_{n-1})M(z,y_{n-1},\eta'')-(z''-x'')\xi''
-z_n^2x_{n-1}\xi_1 \\
& { } &  \ \ \ \ -z_nx_{n-1}^2\xi_1
   - z_{n-1} ^2S_1(z, x_{n-1}, \xi'') -(z_n+x_{n-1})M(z,x_{n-1},\xi'').
\end{eqnarray*}

\par  We will use stationary phase in $z''$ and   $ \eta''$: set
$d_{z''} \tilde{\chi}=0$ and   $ d_{\eta''}\tilde{\chi}=0$, where

\begin{eqnarray*}
d_{z''} \tilde{\chi}&=& \eta''-\xi''+(z_n+y_{n-1})\partial_{z''}
M(z, y_{n-1},
\eta'')-(z_n+x_{n-1})\partial_{z''} M(z, x_{n-1}, \xi'')\\
& { } & \ \ \ \ + z_{n-1}^2 (\partial_{z''}S_1(z, y_{n-1}, \eta'')-
\partial_{z''}S_1(z, x_{n-1}, \xi''))\\
d_{\eta_1} \tilde{\chi} &=& z_1-y_1 + z_n^2y_{n-1} + z_n y_{n-1}^2
+(z_n+y_{n-1})\partial_{\eta_1}M(z, y_{n-1},
\eta'') \\
& { } & \ \ \ \ + z_{n-1}^2 \partial_{\eta_1} S_1(z, y_{n-1}, \eta'')\\
d_{\eta_i} \tilde{\chi} &=& z_i-y_i  +(z_n+y_{n-1})
\partial_{\eta_i}M(z, y_{n-1}, \eta'')+ z_{n-1}^2
\partial_{\eta_i}S_1(z, y_{n-1},
\eta''),  \ \ \ 2 \leq i \leq n-2.
\end{eqnarray*}

\par Notice that $d^2_{z'' \eta''} \tilde{\chi} $ is
nondegenerate. We may solve these equations implicitly for $z''$ and
$ \eta''$ in
terms of the other
variables:

\begin{eqnarray*}
\eta''&=& \xi'' + (z_n+x_{n-1}) \partial_{z''} M(z, x_{n-1},
\xi'')-(z_n+y_{n-1})\partial_{z''}M(z, y_{n-1},\eta'')\\
& { } & \ \ \ \ - z_{n-1}^2 (\partial_{z''}S_1(z, y_{n-1}, \eta'')-
\partial_{z''}S_1(z, x_{n-1},
\xi''))\\
z_i&=&y_i-(z_n+y_{n-1})\partial_{\eta_i}M(z, y_{n-1}, \eta'') -
z_{n-1}^2\partial_{\eta_i}S_1(z, y_{n-1}, \eta'') ,
\ \ \ 2 \leq i \leq n-2;\\
z_1&=&y_1-z_n^2y_{n-1} -z_ny_{n-1}^2-(z_n+y_{n-1})
\partial_{z_1}M(z, y_{n-1},
\eta'')- z_{n-1}^2 \partial_{\eta_1}S_1(z, y_{n-1}, \eta'').
\end{eqnarray*}

\par We have that $\partial_{z''}S_1(z, y_{n-1},
\eta'')- \partial_{z''}S_1(z, x_{n-1}, \xi'')$ vanishes at
$y_{n-1}= x_{n-1}$ and $\eta''=\xi''$, so we can write

\begin{eqnarray*}
\partial_{z''}S_1(z, y_{n-1}, \eta'')- \partial_{z''}S_1(z,
x_{n-1}, \xi'')= (y_{n-1}- x_{n-1}) \partial_{z'' y_{n-1}}S_1(z,
y_{n-1}, \eta'') \\ + (\eta''-\xi'') \partial_{z''
\eta''}S_1(z,y_{n-1}, \eta'')
\end{eqnarray*}

\par In a similar way,

\begin{eqnarray*}
  (z_n&+&y_{n-1})  \partial_{z''} M(z, y_{n-1},
\eta'')-(z_n+x_{n-1})\partial_{z''}M(z, x_{n-1},\xi'') \\& = &
(y_{n-1}-x_{n-1})[ \partial_{z''}M(z, y_{n-1}, \eta'') +
(z_n+y_{n-1})\partial_{z'' y_{n-1}}M(z, y_{n-1}, \eta'')] \\ & &+
(\eta''- \xi'')(z_n+ y_{n-1})
\partial_{z'' \eta''}M(z,y_{n-1}, \eta )
\end{eqnarray*}

\par Thus, $\eta''- \xi''$ becomes

\begin{eqnarray*}
  \eta''-\xi''&=&-(y_{n-1}-x_{n-1})[
\partial_{z''}M(z, y_{n-1}, \eta'') \\
& &+ (z_n+y_{n-1})\partial_{z''
y_{n-1}}M(z, y_{n-1}, \eta'')  + z_{n-1}^2 \partial_{z''
y_{n-1}}S_1(z, y_{n-1}, \eta'')] \\
&\times&[I + (z_n+y_{n-1})
\partial_{z'' \eta''}M(z, y_{n-1}, \eta'')+ z_{n-1}^2 \partial_{z'' 
y_{n-1}}S_1(z, y_{n-1}, \eta'')]^{-1}
\end{eqnarray*}

\par The phase $\tilde {\chi}$ then becomes:

\begin{eqnarray*}
\tilde{\tilde
{\chi}}&=&(x''-y'')\cdot\xi''+z_n^2(y_{n-1}-x_{n-1})\xi_1+
z_n(y_{n-1}^2-x_{n-1}^2)\xi_1\\
& & +(y_{n-1} -x_{n-1})[z_{n-1}^2
\partial_{y_{n-1}}S_1(\cdot)  + M(\cdot)
+ (z_n + y_{n-1}) \partial_{y_{n-1}}M(\cdot)]
\\
&=&(x''-y'')\cdot\xi'' + (y_{n-1}-x_{n-1}) \{ \left(z_n^2 + z_n
(y_{n-1}+x_{n-1})\right)\xi_1\\
& &+ z_{n-1}^2 (\partial_{y_{n-1}}S_1(\cdot) + S_5(\cdot))\\
& &  + (z_n +
y_{n-1})[z_{n-1}S_6(\cdot) + (z_n + y_{n-1})S_7(\cdot) +
\partial_{y_{n-1}}M(\cdot) ] \}
\end{eqnarray*}
where $(\cdot)= (y'', z_{n-1}, z_n, y_{n-1}, \xi'')$  and the
amplitude becomes $\tilde{a}\in S^{2m+1}$.\\

\par Let
\begin{eqnarray*}
\ N (\cdot)&=&z_{n-1}^2 (\partial_{y_{n-1}}S_1(\cdot) + S_5(\cdot))\\
\ \ & &+ (z_n + y_{n-1})[z_{n-1}S_6(\cdot) + (z_n + y_{n-1})S_7(\cdot)
+\partial_{y_{n-1}}M(\cdot)]
\end{eqnarray*}
and

\begin{equation*}
\ Q(\cdot)=z_{n-1}S_6(\cdot) + (z_n + y_{n-1})S_7(\cdot) +
\partial_{y_{n-1}}M(\cdot)
\end{equation*}
with
\begin{eqnarray}\label{old-40}
\partial_{y_{n-1}}M&=&z^2_{n-1}\partial_{y_{n-1}}S_5+ z_{n-1}S_6+
z_{n-1}(z_n+y_{n-1})\partial_{y_{n-1}}S_6\nonumber\\
& &+
2(z_n+y_{n-1})S_7+(z_n+y_{n-1})^2\partial_{y_{n-1}}S_7.
\end{eqnarray}

\par Repeating the argument from \S5,
we make a  singular change of variables,

\begin{eqnarray*}
\theta_i&=&\xi_i,  \ \ \  1 \leq i \leq n-2,\\
\theta_{n-1}&=& -(z_n^2\xi_1 + z_n (y_{n-1}+x_{n-1})\xi_1 +
N(\cdot)).
\end{eqnarray*}

\par We have
\begin{eqnarray*}
\nabla_{z_{n-1}, z_n} \theta_{n-1}&=&-\big(2z_{n-1} [
\partial_{y_{n-1}}S_1(\cdot)
+ S_5(\cdot)]+ (z_n + y_{n-1}) \partial_{z_{n-1}}Q(\cdot),\\
&{ }&\ \ \ \  2z_n \theta_1 +(x_{n-1}+y_{n-1})\theta_1+
\partial_{z_n}N(\cdot) \big),
\end{eqnarray*}
so that $|\nabla_z \theta_{n-1}|=0 $ \ \ {\rm iff }
\begin{equation*}
z_{n-1}=- \frac{(z_n+y_{n-1})\partial_{z_{n-1}}Q(\cdot)}{2 (
\partial_{y_{n-1}}S_1(\cdot)
+ S_5(\cdot))}
\end{equation*}
and
\begin{equation*}
z_n=-\frac{x_{n-1}+y_{n-1}}{2}-\frac{\partial_{z_n}N(\cdot)}{2
\theta_1}.
\end{equation*}

At these points,

\begin{eqnarray*}
\theta_{n-1}&=& \frac{(x_{n-1}+y_{n-1})^2}{4} \theta_1- \frac{1}{4
\theta_1 }(\partial_{z_n}N(\cdot))^2\\
& { } & - (z_n+y_{n-1})^2\frac{(\partial_{z_{n-1}}Q)^2}{4
(\partial_{y_{n-1}}S_1(\cdot) +
S_5(\cdot))} - (z_n+y_{n-1}) Q(\cdot)\\
&:=& \frac{(x_{n-1}+y_{n-1})^2}{4} \theta_1 +P(y'', y_{n-1},
x_{n-1}, \theta'').
\end{eqnarray*}
where $P(y'', y_{n-1}, x_{n-1}, \theta'')$ is

\begin{eqnarray}\label{old-43}
\tilde{P}(\cdot):= &-& \frac{1}{4 \theta_1
}(\partial_{z_n}N(\cdot))^2\\ &-&
(z_n+y_{n-1})^2\frac{(\partial_{z_{n-1}}Q)^2}{4
(\partial_{y_{n-1}}S_1(\cdot) + S_5(\cdot))} - (z_n+y_{n-1})
Q(\cdot)\nonumber
\end{eqnarray}
pushed forward under the previous change of variables.
\medskip

\par Letting $\so^{C_0}(z,y)$ denote the critical set in the $z,y$
variables, i.e.,  for
$C\subset T^*Z\times T^*Y$, note that
$\tilde{P}(\cdot)|_{\so^{C_0}(z,y)}=0$ and $\nabla
\tilde{P}(\cdot)|_{\so^{C_0}(z,y)}=0$, so that
\begin{equation}\label{fgp-1}
P|_{ \{
x_{n-1}=y_{n-1} \}}=0 \hbox{ and }\nabla P|_{ \{ x_{n-1}=y_{n-1} \}
}=0.
\end{equation}
\medskip

\par As in the model case, it follows that $K_{F^*F}(x,y)$ has an
oscillatory integral representation,

\[\int_{\R^n} e^{i\{(x''-y'')\cdot\theta'' +
(x_{n-1}-y_{n-1})\theta_{n-1} +\frac{\rho}{\theta_1}(\theta_{n-1}-
\frac{(x_{n-1}+y_{n-1})^2}{4} \theta_1 -P)\}} a(x,y, \theta, \rho)
d\theta d \rho,\]
where $a \in S^{2m,-1}(2n-2,n-1,1)$.  On the region $\{|\rho|\le c
|\xi|\}$, the new
phase function, $\psi(x,y;\theta;\rho)$ is a multi-phase for  a pair
$(\Delta', {\tilde C}')$ in the sense of Def.\ref{old-def3.6}:
$\psi(x,y;\theta;0)$
parametrizes the diagonal Lagrangian $\Delta'$ and $\psi(x,y;\theta;\rho)$
parametrizes a Lagrangian ${\tilde C}'$.  Hence, the contribution to
$F^*F$ from this
region is in
$I^{p,l} (\Delta,\tilde C)$ for some
$p,l\in \R$, and the orders of $F^*F$  are computed using Remark 3.7
in the same way
as for $K^0$ in the model case in \S5, so that $p=2m-\frac12$ and
$l=\frac12$. On the
other hand, the contribution from $\{|\rho|\ge c |\xi|\}$ is handled
in the same way
as for
$K^1$  for the model case in \S5, giving an element of $I^{2m-\frac12}(\tilde
C)\subset I^{2m-\frac12,\frac12}(\Delta,\tilde C)$.

\par Next, we show that ${\tilde C}$ is a two-sided fold; the fact
that $\tilde C$ is
symmetric just follows  from $\Delta\cup\tilde C= C^t\circ C$. We have:

\begin{eqnarray*}
\tilde C = \Big\{ \big(&x''&, x_{n-1}, \theta'', \theta_{n-1} -
\frac{x_{n-1}+y_{n-1}}{2}  \rho -
\frac{\rho}{\theta_1} \partial_{x_{n-1}} P;\\
   &y''&, y_{n-1},
\theta'', \theta_{n-1} + \frac{x_{n-1}+y_{n-1}}{2} \rho +
\frac{\rho}{\theta_1}
\partial_{y_{n-1}} P\big):\\
& &
\theta_{n-1}- \frac{(x_{n-1} + y_{n-1})^2}{4}\theta_{1} - P=0,
x_{n-1}-y_{n-1} + \frac{\rho}{\theta_1} =0,\\
& & x_i-y_i + \frac{\rho}{\theta_1} \partial_{\theta_i} P=0,\ \   2
\leq i \leq n-2,\\
& &  x_1-y_1 -\frac{\rho}{\theta_1^2}(\theta_{n-1}-
\frac{(x_{n-1}+y_{n-1})^2}{4} \theta_1 -P) \\
& &- \frac{\rho}{\theta_1}
\frac{(x_{n-1}+y_{n-1})^2}{4} -\frac{\rho}{\theta_1}
\partial_{\theta_1}P=0 \Big\}.
\end{eqnarray*}

\par The coordinates on $\tilde C$ are $(x'', x_{n-1}, y_{n-1},\theta_1,
\theta''')$. Note that
$\rho=-(x_{n-1}-\nolinebreak y_{n-1})\theta_1$ and $\theta_{n-1}=
\frac{(x_{n-1} + y_{n-1})^2}{4}\theta_{1} + P$; thus, $y_i$ and
$y_1$ become:

$$y_i= x_i-(x_{n-1}-y_{n-1}) \partial_{\theta_i}P,  \ \  2 \leq i
\leq n-2,$$
and
$$y_1=x_1 + \frac{(x_{n-1}-y_{n-1})(x_{n-1}+y_{n-1})^2}{4} +
(x_{n-1}-y_{n-1}) \partial_{\theta_1}P.$$

\par We now consider the projections $\pi_R$ and $\pi_L$.
We  have

\begin{eqnarray*}
\pi_R \big(x_1, x''', y_{n-1}, &\theta_1&, \theta''', x_{n-1}\big)\\
&=&\Big(x_1+\frac{(x_{n-1}-y_{n-1})(x_{n-1}+y_{n-1})^2}{4} +
(x_{n-1}+y_{n-1})\partial_{\theta_1}P,\\
   &\ \ &\ \ x''' -(x_{n-1}-y_{n-1})
\partial_{\theta'''} P, y_{n-1};\theta_1, \theta''', \\
&\ \ &
\frac{(x_{n-1} + y_{n-1})(3y_{n-1}-x_{n-1})}{4}\theta_{1} + P -
(x_{n-1}-y_{n-1})\partial_{y_{n-1}} P\Big)
\end{eqnarray*}
and

\begin{eqnarray*}
\pi_L\big( &x''& ,\ \  x_{n-1}\ \  ,\ \  \theta_1\ \ ,\ \  \theta'''\
\ ,
\ \  y_{n-1}\ \
\big)\\ &=& \left(x'', x_{n-1};
\theta_1,
\theta''',
\ \  \frac{(x_{n-1} +
y_{n-1})(3x_{n-1}-y_{n-1})}{4}\theta_{1} + P  + (x_{n-1}-y_{n-1})
\partial_{x_{n-1}} P \right).
\end{eqnarray*}

\par  One easily computes

\[{d \pi_L}= \left(\begin{array}{ccccc}
I_{n-1} & 0 & 0 & 0 & 0 \\
0 & 1 & 0 & 0 & 0\\
0 & 0 & 1 & 0  & 0\\
0 & 0 &0 & I_{n-3} & 0  \\
0 & \cdot & \cdot & \cdot &  D
\end{array} \right)\]
where $D= \frac{x_{n-1}-y_{n-1}}{2} \theta_1 +
\partial_{y_{n-1}}P-\partial_{x_{n-1}}P +
(x_{n-1}-y_{n-1})\partial^2_{x_{n-1}y_{n-1}}P$ and 
Ker\nolinebreak$\quad d \pi_L$
is spanned by $ \frac{\partial}{\partial y_{n-1}}$.

\par Similarly,

\[{d \pi_R}= \left(\begin{array}{cccccc}
1 & \cdot & \cdot & \cdot & \cdot & \cdot \\
A & B & \cdot & \cdot & \cdot & \cdot \\
0 &  0  & 1 & 0 & 0  & 0 \\
0 & 0 &  0 & 1 & 0 & 0  \\
0 & 0 & 0 & 0 & I_{n-3} & 0 \\
\cdot &  \cdot & \cdot & \cdot &  \cdot & -D
\end{array} \right)\]
where $A=-(x_{n-1}-y_{n-1})\partial_{\theta''' x_1}^2P$ and
$B=I_{n-3}-(x_{n-1}-y_{n-1}) \partial_{\theta ''' x'''}^2P$.
Hence,  det $d \pi_R= D$, as well,  and  $\Ker d\pi_R$ is spanned
by $\frac{\partial\quad}{\partial x_{n-1}}$.

\par We have:
$$\partial_{y_{n-1}}D= -\frac{1}{2} \theta_1 +
\partial^2_{y^2_{n-1}}P - 2 \partial^2_{x_{n-1}y_{n-1}}P
  + (x_{n-1}-y_{n-1})\partial^3_{x_{n-1}y^2_{n-1}}P$$
and $$\partial_{x_{n-1}}D= \frac{1}{2} \theta_1 -
\partial^2_{x^2_{n-1}}P + 2 \partial^2_{x_{n-1}y_{n-1}}P
  + (x_{n-1}-y_{n-1})\partial^3_{x^2_{n-1}y_{n-1}}P$$

Thus, to show that $\pi_R$ and $\pi_L$ are both folds, with
$\so^{\tilde{C}}=\{x_{n-1}=y_{n-1}\}=\Delta\cap\tilde{C}$, it
suffices to prove that
\begin{equation}\label{fgp-2}
|\partial^2_{y^2_{n-1}}P|+|
\partial^2_{x_{n-1}y_{n-1}}P|+|\partial^2_{x^2_{n-1}}P|\le\frac18
|\theta_1|
\end{equation}
on a neighborhood of $\tilde{c}_0$;
to do this, we need  to show that
$\partial^2_{y^2_{n-1}}\tilde{P},
\partial^2_{x_{n-1}y_{n-1}}\tilde{P}$
and $\partial^2_{x^2_{n-1}}\tilde{P}$ satisfy the
same estimate near the corresponding point.

\par Let $T=\frac{(\partial_{z_{n-1}}Q)^2}{4
(\partial_{y_{n-1}}S_1(\cdot) + S_5(\cdot))} $.
 From (\ref{old-43}), we see that

\begin{eqnarray*}
\partial_{y_{n-1}} \tilde{P}&=& -\frac{1}{2\theta_1} \partial_{z_n}N
  \partial_{z_n y_{n-1}}N -2(z_n+y_{n-1})(\partial_{y_{n-1}}z_n
  +1)T \\
  & &- (z_n+y_{n-1})^2
\partial_{y_{n-1}} T-(z_n+y_{n-1})\partial_{y_{n-1}}Q- 
(\partial_{y_{n-1}}z_n+1)Q
\end{eqnarray*}
and

\begin{eqnarray*}
\partial_{x_{n-1}} \tilde{P}&=& -\frac{1}{2\theta_1} \partial_{z_n}N
  \partial_{z_n x_{n-1}}N -2(z_n+y_{n-1})(\partial_{x_{n-1}}z_n
  )T\\
  & &- (z_n+y_{n-1})^2
\partial_{x_{n-1}} T-(z_n+y_{n-1})\partial_{x_{n-1}}Q- 
(\partial_{x_{n-1}}z_n)Q.
\end{eqnarray*}
Thus
\begin{eqnarray*}
   \partial^2_{y^2_{n-1}}\tilde{P}& =&
-\frac{1}{2\theta_1}(\partial^2_{z_ny_{n-1}}N)^2-\frac{1}{2\theta_1}\partial_{z_n}N
\partial^3_{z_ny^2_{n-1}}N \\
& &-2 (\partial_{y_{n-1}}z_n+1)^2 T-2(z_n+y_{n-1})(\partial^2_{y^2_{n-1}}z_n
)T\\
& &- 4(z_n+y_{n-1})(\partial_{y_{n-1}}z_n
  +1)\partial_{y_{n-1}}T -(z_n+y_{n-1})^2
\partial^2_{y^2_{n-1}} T\\
& & - 2(\partial_{y_{n-1}}z_n+1)\partial_{y_{n-1}}Q
-(z_n+y_{n-1})\partial^2_{y^2_{n-1}}Q -(\partial^2_{y^2_{n-1}}z_n)Q.
\end{eqnarray*}
  Similarly,

  \begin{eqnarray*}
  \partial^2_{x^2_{n-1}}\tilde{P} &=&
-\frac{1}{2\theta_1}(\partial^2_{z_nx_{n-1}}N)^2-\frac{1}{2\theta_1}\partial_{z_n}N
\partial^3_{z_nx^2_{n-1}}N \\
& &-2 (\partial_{x_{n-1}}z_n)^2 T-2(z_n+y_{n-1})(\partial^2_{x^2_{n-1}}z_n
)T\\
& &- 4(z_n+y_{n-1})(\partial_{x_{n-1}}z_n
  )\partial_{x_{n-1}}T -(z_n+y_{n-1})^2
\partial^2_{x^2_{n-1}} T\\
& & - 2(\partial_{x_{n-1}}z_n)\partial_{x_{n-1}}Q
-(z_n+y_{n-1})\partial^2_{x^2_{n-1}}Q -(\partial^2_{x^2_{n-1}}z_n)Q
\end{eqnarray*}
and
\begin{eqnarray*}
  \partial^2_{y_{n-1}x_{n-1}}\tilde{P} &=&
-\frac{1}{2\theta_1}(\partial^2_{z_ny_{n-1}}N)
\partial^2_{z_nx_{n-1}}N-\frac{1}{2\theta_1}\partial_{z_n}N
\partial^3_{z_ny_{n-1}x_{n-1}}N \\
& &-2 (\partial_{x_{n-1}}z_n)(\partial_{y_{n-1}}z_n+1) T-
2(z_n+y_{n-1})(\partial^2_{x_{n-1}y_{n-1}}z_n )T\\
& &-
2(z_n+y_{n-1})(\partial_{y_{n-1}}z_n
  +1)\partial_{x_{n-1}}T - 2(z_n+y_{n-1})(\partial_{x_{n-1}}z_n
  )\partial_{y_{n-1}}T\\
& &-(z_n+y_{n-1})^2
\partial^2_{y_{n-1}x_{n-1}} T - (\partial_{x_{n-1}}z_n)
\partial_{y_{n-1}}Q\\
& &-(z_n+y_{n-1})\partial^2_{y_{n-1}x_{n-1}}Q
-(\partial^2_{y_{n-1}x_{n-1}}z_n) Q -
(\partial_{y_{n-1}}z_n+1)\partial_{x_{n-1}}Q.
\end{eqnarray*}

\par We have
\begin{eqnarray*}
\partial_{z_n}N&=&z_{n-1}^2(\partial^2_{z_ny_{n-1}}S_1 + \partial_{z_n}S_5)\\
& &+ (z_n+y_{n-1})
[z_{n-1}\partial_{z_n}S_6 + S_7 + (z_n+y_{n-1})\partial_{z_n}S_7 +
\partial^2_{y_{n-1}z_n}M] + Q.
\end{eqnarray*}
from which it follows that
\begin{eqnarray}\label{old-49}
\partial^2_{z_n y_{n-1}}N&=& 2 
z_{n-1}\partial_{y_{n-1}}z_{n-1}(\partial^2_{z_ny_{n-1}}S_1 +
\partial_{z_n}S_5)+ z^2_{n-1}\partial_{y_{n-1}}(\partial^2_{z_ny_{n-1}}S_1 +
\partial_{z_n}S_5)\nonumber\\
& &+ (\partial_{y_{n-1}}z_n+1)[z_{n-1}\partial_{z_n}S_6 + S_7 + 
(z_n+y_{n-1})\partial_{z_n}S_7 +
\partial^2_{y_{n-1}z_n}M] \nonumber\\
& &+ (z_n+y_{n-1})\partial_{y_{n-1}}[z_{n-1}\partial_{z_n}S_6 + S_7 + 
(z_n+y_{n-1})\partial_{z_n}S_7 +
\partial^2_{y_{n-1}z_n}M]\nonumber\\
& & +\partial_{y_{n-1}}Q
\end{eqnarray}
and
\begin{eqnarray}\label{old-50}
  \partial^2_{z_n x_{n-1}}N&=& 2 
z_{n-1}\partial_{x_{n-1}}z_{n-1}(\partial^2_{z_ny_{n-1}}S_1 +
\partial_{z_n}S_5)+ z^2_{n-1}\partial_{x_{n-1}}(\partial^2_{z_ny_{n-1}}S_1 +
\partial_{z_n}S_5)\nonumber\\
& &+ (\partial_{x_{n-1}}z_n)[z_{n-1}\partial_{z_n}S_6 + S_7 + 
(z_n+y_{n-1})\partial_{z_n}S_7 +
\partial^2_{y_{n-1}z_n}M] \nonumber\\
& &+ (z_n+y_{n-1})\partial_{x_{n-1}}[z_{n-1}\partial_{z_n}S_6 + S_7\nonumber\\
& & + (z_n+y_{n-1})\partial_{z_n}S_7 +
\partial^2_{y_{n-1}z_n}M] +\partial_{x_{n-1}}Q.
\end{eqnarray}
Similarly,
\begin{eqnarray}\label{old-51}
\partial_{y_{n-1}}Q&=&\partial_{y_{n-1}}z_{n-1}S_6+
z_{n-1}\partial_{y_{n-1}}S_6+
(\partial_{y_{n-1}}z_n+1)S_7\nonumber\\
& &+(z_n+y_{n-1})\partial_{y_{n-1}}S_7+\partial^2_{y^2_{n-1}}M
\end{eqnarray}
and
\begin{eqnarray}\label{old-52}
\partial_{x_{n-1}}Q&=&\partial_{x_{n-1}}z_{n-1}S_6+
z_{n-1}\partial_{x_{n-1}}S_6\nonumber\\
& &+
(\partial_{x_{n-1}}z_n)S_7+(z_n+y_{n-1})\partial_{x_{n-1}}S_7+\partial^2_{y_{n-1}x_{n-1}}M.
\end{eqnarray}
\vskip.3in

Using (\ref{old-40}), we have
\begin{eqnarray}\label{old-53}
\partial^2_{y^2_{n-1}}M&=&2z_{n-1}(\partial_{y_{n-1}}z_{n-1})
\partial_{y_{n-1}}S_5+z^2_{n-1}\partial^2_{y^2_{n-1}}S_5+
\partial_{y_{n-1}}z_{n-1}S_6 \nonumber\\
& &+
z_{n-1}\partial_{y_{n-1}}S_6+\partial_{y_{n-1}}z_{n-1}(z_n+y_{n-1})
\partial_{y_{n-1}}S_6 \nonumber\\
& &+z_{n-1}(\partial_{y_{n-1}}z_n+1)
\partial_{y_{n-1}}S_6+z_{n-1}(z_n+y_{n-1})\partial^2_{y^2_{n-1}}S_6
\nonumber\\
& &+
2(\partial_{y_{n-1}}z_n+1)S_7+2(z_n+y_{n-1})\partial_{y_{n-1}}S_7\\
& &
+2(z_n+y_{n-1})(\partial_{y_{n-1}z_n}+1)\partial_{y_{n-1}}S_7
+(z_n+y_{n-1})^2\partial^2_{y^2_{n-1}}S_7 \nonumber
\end{eqnarray}
and
\begin{eqnarray}\label{old-54}
\partial^2_{y_{n-1}x_{n-1}}M&=&2z_{n-1}(\partial_{x_{n-1}}z_{n-1})
\partial_{y_{n-1}}S_5+z^2_{n-1}\partial^2_{y_{n-1}x_{n-1}}S_5+
\partial_{x_{n-1}}z_{n-1}S_6 \nonumber\\
& &+
z_{n-1}\partial_{x_{n-1}}S_6+\partial_{x_{n-1}}z_{n-1}(z_n+y_{n-1})
\partial_{y_{n-1}}S_6
\nonumber\\
& &+z_{n-1}(\partial_{x_{n-1}}z_n)
\partial_{y_{n-1}}S_6+z_{n-1}(z_n+y_{n-1})\partial^2_{y_{n-1}x_{n-1}}S_6
\nonumber\\
& &+2(\partial_{x_{n-1}}z_n)S_7+2(z_n+y_{n-1})\partial_{x_{n-1}}S_7
\\
& &+2(z_n+y_{n-1})(\partial_{x_{n-1}}z_n)\partial_{y_{n-1}}S_7
+(z_n+y_{n-1})^2\partial^2_{y_{n-1}x_{n-1}}S_7.\nonumber
\end{eqnarray}
Similarly,
\begin{eqnarray}\label{old-55}
\partial^2_{y_{n-1}z_n}M&=&z_{n-1}^2\partial_{y_{n-1}z_n}S_5+z_{n-1}\partial_{z_n}S_6 
+
z_{n-1}\partial^2_{y_{n-1}}S_6
  \nonumber\\
& &+z_{n-1}(z_n+y_{n-1})\partial^2_{y_{n-1}z_n}S_6+
2S_7+2(z_n+y_{n-1})\partial_{z
_n}S_7
\nonumber\\
& 
&+2(z_n+y_{n-1})\partial_{y_{n-1}}S_7+(z_n+y_{n-1})^2\partial^2_{y_{n-1}z_n}S_7.
\end{eqnarray}

\par Using (\ref{old-49}) - (\ref{old-55})  and  the fact that 
$\partial_{z_n}N, Q, \partial_{y_{n-1}}M,
S_6$ and $S_7$ are all small near $c_0$, we obtain that
$\partial^2_{y_{n-1}x_{n-1}}M, \ \  \partial^2_{y^2_{n-1}}M,
\ \ \partial^2_{y_{n-1}z_n}M, \ \ \partial_{x_{n-1}}Q$,
$\partial_{y_{n-1}}Q, \ \ \partial^2_{z_ny_{n-1}}N$ and
$ \ \ \ \ \ \ \partial^2_{z_nx_{n-1}}N$ are as well.  Thus,
$\partial^2_{y^2_{n-1}}\tilde{P},
\partial^2_{x^2_{n-1}}\tilde{P}$ and
$\partial^2_{x_{n-1}y_{n-1}}\tilde{P}$ are also small, from which
it follows that (\ref{fgp-2}) holds, and thus  $\tilde{C}$ is a
symmetric, two-sided fold, with
$\Delta\cap\tilde{C}={\cal S}_1^{\tilde C}$.

\begin{remark}
{\rm An interesting question, particularly relevant for the marine
seismic imaging
problem, is whether (under an ellipticity assumption on $F$), } $F^*F$
{\rm can be inverted, at least modulo operators mapping} $H^s
\rightarrow H^{s-2m + \delta}$ {\rm for some} $\delta >0$. {\rm
Note that, by Remark \ref{old-rem3.4}}, {\rm the order of} $F^*F$ {\rm on}
$\tilde{C} \setminus \Delta$ {\rm is } $\frac{1}{2}$ {\rm
less than on} $\Delta \setminus \tilde{C} $. {\rm  However, the
standard technique of parabolic decomposition (see \cite{gruh2})
only results in a decomposition} $F^*F= T_1 +T_2$ {\rm with} $T_1
\in I_{\frac{1}{2}, \frac{1}{2}}^{2m} (\Delta)$ {\rm and} $T_2 \in
I_{\frac{1}{2}, \frac{1}{2}}^{2m} (\tilde{C})$ {\rm .
Since} $\tilde{C}$ {\rm is a two-sided fold, there is further loss
of} $\frac{1}{6}$ {\rm derivatives in terms of Sobolev mapping
properties  \cite{meta}.  We will return to this question in  \cite{fgp}.}

\end{remark}

\vskip.2in

\noindent{\sc Department of Mathematics and Statistics}

\noindent{\sc Rochester  Institute of Technology

\noindent{\sc Rochester, NY 14623}

\noindent{\tt{rxfsma@rit.edu}}

\vskip.2in

\noindent{\sc Department of Mathematics

\noindent{\sc University of Rochester

\noindent{\sc Rochester, NY 14627}

\noindent{\tt{allan@math.rochester.edu}}

\end{document}